\documentclass[journal,twoside,web]{ieeecolor}
\usepackage{generic}

\usepackage{amsmath,amssymb,amsfonts,mathrsfs,mathtools}
\usepackage{algorithm,algorithmic}
\usepackage{graphicx}
\usepackage{changepage}
\usepackage{cite, comment}

\makeatletter
\let\NAT@parse\undefined
\makeatother
\usepackage[
    colorlinks=false,
    pdfborder={1 1 1},
    linkbordercolor={1 0 0},
    citebordercolor={1 0 0},
]{hyperref}

\usepackage{textcomp}
\def\BibTeX{{\rm B\kern-.05em{\sc i\kern-.025em b}\kern-.08em
    T\kern-.1667em\lower.7ex\hbox{E}\kern-.125emX}}
\DeclareMathOperator*{\argmin}{arg\,min}
\DeclareMathOperator*{\rank}{\mathrm{rank}}

\DeclareMathOperator*{\blkdiag}{\mathrm{blkdiag}}
\DeclareMathOperator*{\sym}{\mathrm{sym}}

\DeclarePairedDelimiter{\norm}{\lVert}{\rVert}
\newtheorem{theorem}{Theorem}
\newtheorem{lemma}{Lemma}
\newtheorem{assumption}{Assumption}
\newtheorem{remark}{Remark}
\newtheorem{example}{Example}
\newtheorem{corollary}{Corollary}
\newtheorem{definition}{Definition}

\usepackage{tikz}
\usetikzlibrary{matrix,positioning}

\begin{document}
\title{Data-driven stabilization of nonlinear systems via descriptor embedding}
\author{
	Mohammad Alsalti, Claudio De Persis, Victor G. Lopez, and Matthias A. M\"uller
	\thanks{This work has received funding from the European Research Council (ERC) under the European Union’s Horizon 2020 research and innovation programme (grant agreement No 948679).}
    \thanks{M. Alsalti thanks the Graduiertenakademie at Leibniz University Hannover for funding his research stay the University of Groningen.}
	\thanks{M. Alsalti, V. G. Lopez and M. A. M\"uller are with the Leibniz University Hannover, Institute of Automatic Control, 30167 Hannover, Germany (email: \{alsalti,lopez,mueller\}@irt.uni-hannover.de).}
	\thanks{C. De Persis is with the ENTEG and the J. C. Willems Center for Systems and Control, University of Groningen, 9747 AG Groningen, The Netherlands (e-mail: c.de.persis@rug.nl).}
}

\maketitle

\begin{abstract}
	We introduce the notion of descriptor embedding for nonlinear systems and use it for the data-driven design of stabilizing controllers. Specifically, we provide sufficient data-dependent LMI conditions which, if feasible, return a stabilizing nonlinear controller of the form $\mathnormal{u=K(x)Z(x)}$ where $\mathnormal{K(x)}$ belongs to a polytope and $\mathnormal{Z}$ is a user-defined function. The proposed method is then extended to account for the presence of uncertainties and noisy data. Furthermore, a method to estimate the resulting region of attraction is given using only data. Simulation examples are used to illustrate the results and compare them to existing methods from the literature.
\end{abstract}

\begin{IEEEkeywords}
	Control design, data-driven control, linear matrix inequalities, nonlinear systems
\end{IEEEkeywords}

\section{Introduction}
\label{sec:intro}
\IEEEPARstart{D}{}irect data-driven analysis and control \cite{markovsky2021behavioral,florian23ddcontrol} refers to the use of historical input-state (or input-output) data collected from a system to directly study its properties or obtain stabilizing controllers, thus avoiding explicit model identification steps. The origins of many of the recent works on this topic can be traced back to results from the behavioral approach to systems theory~\cite{willems1997book}. For instance, for a linear time-invariant (LTI) system, a so-called \textit{data-based representation} of all its trajectories is given by the image of a Hankel matrix of data that satisfies certain richness conditions \cite{willems2005note,markovsky2022identifiability,Alsalti2025DBrep}. This was later used for the direct data-driven design of (robust and/or optimal) stabilizing state-feedback controllers \cite{depersis2019formulas,vanwaarde2023informativity,bisoffi2022data,van2020noisy}, predictive controllers \cite{berberich2024overview} and output-feedback controllers~\cite{alsalti2023notes,Li2024} for LTI systems.

Data-based representations for classes of nonlinear systems can also be obtained by exploiting system theoretic properties and/or by assuming the availability of user-defined basis functions that span the dynamics. For instance, Hammerstein-Wiener systems were studied in \cite{berberich2020trajectory}, Volterra systems in \cite{rueda2020data}, linear parameter-varying systems in \cite{verhoek2021fundamental}, feedback linearizable systems in \cite{alsalti2023data} and generalized bilinear systems in \cite{markovsky2023bilinear}. Direct design of stabilizing controllers for classes of nonlinear systems has also received significant attention (see \cite{martin2023guarantees,depersis2023annrev} and the references therein). Of those techniques, many mainly relied on the use of basis functions, e.g., for stabilization of polynomial systems using sum-of-squares programming as in \cite{Guo2022polynomial}, or to cancel/dominate the nonlinearities as in \cite{DePersis2023}, or to enforce contraction \cite{hu2024enforcing}. Alternatively, one may use a lifting technique and study the dynamics in a linear or bilinear state-space using, e.g., the Koopman operator \cite{mauroy2020,robin2024koopman}. However, the lifted system is, in general, infinite dimensional and the so-called Koopman eigenfunctions are usually not determined systematically but instead are approximated from data~\cite{Haseli2022}.

Methods that rely on (approximate) cancellation and/or domination of the system's nonlinearities can potentially perform poorly in the presence of uncertainties and disturbances. Instead, it would be advantageous to exploit the inherent nonlinearities in the system in order to achieve better performance and/or robustness properties. To this end, in this work we introduce a notion of \textit{descriptor embedding for nonlinear systems}. This allows for systematic analysis and control design for nonlinear systems using tools originally developed for (linear parameter-varying) descriptor systems, see, e.g.,~\cite{dai1989}. Considering a descriptor viewpoint allows for obtaining controllers of the form $u=K(x)Z(x)$, which is different from the non-descriptor linear parameter-varying embedding, e.g., in \cite{verhoek2023general}, that returns (less expressive) controllers of the form $u=K(p)x$ (where $p$ is a scheduling parameter that must be both well-defined and available for measurement). In contrast, our technique directly uses the state as the measured variable. Our proposed embedding technique is applicable in model-based and model-free fashion and the goal of this work is to highlight its usefulness in the model-free domain by designing stabilizing controllers for nonlinear systems purely from data.

The following are the main contributions of this paper. First, we introduce the notion of descriptor embedding for nonlinear systems. This is a new lifting procedure wherein a stabilizing controller is designed for a higher dimensional descriptor system and then shown to be stabilizing for the original nonlinear system. The resulting controller takes the form $u=K(x)Z(x)$ where $K(x)$ belongs to a polytope whose vertices are found by solving a set of data-dependent linear matrix inequalities (LMIs). Notably, the proposed approach allows for obtaining new (potentially global) stabilization techniques for a large class of nonlinear systems which are \emph{not} based on nonlinearity cancellation. Second, the proposed approach is extended to account for uncertainties (e.g., neglected nonlinearities) and noisy data. Finally, we illustrate the proposed method using simulation examples and show that it often outperforms existing data-driven design techniques that are based on nonlinearity cancellation/minimization. This is due to the controller structure which allows for a state-dependent control gain and, hence, results in better closed-loop performance.

The remainder of this paper is structured as follows: Section~\ref{sec:prelim} introduces the notation and useful preliminaries. In Section~\ref{sec:embedding} we introduce the notion of descriptor embedding for nonlinear systems and later, in Section~\ref{sec:nominal}, use this notion to develop data-dependent nonlinear control design techniques for stabilizing the origin of nonlinear systems. In Section~\ref{sec:extensions}, we extend the results of the previous sections to accommodate the presence of noise in the data and neglected nonlinearities. Finally, Section~\ref{sec:conc} concludes the paper.
\section{Notation and preliminaries}\label{sec:prelim}
The sets of natural and real numbers are denoted by $\mathbb{N},\mathbb{R}$, respectively. We use $I_m$ to denote an $m\times m$ identity matrix and $0_{n\times m}$ to denote an $n\times m$ matrix of zeros. We omit the subscript when the dimensions are clear from the context. A symmetric positive definite (semi-definite) matrix is denoted $P\succ0\,(\succeq0)$. Similarly, negative definite (semi-definite) matrices are denoted by $P\prec0\,(\preceq0)$. For a square matrix $M$, we define $\sym(M)\coloneqq M+M^\top$. Blocks of symmetric matrices are denoted by $\star$ whenever they are inferred from context, e.g.,
\[
\begin{bmatrix}
	M_1 & M_2^\top \\ M_2 & M_3
\end{bmatrix} = \begin{bmatrix}
	M_1 & M_2^\top\\ \star & M_3
\end{bmatrix} = \begin{bmatrix}
	M_1 & \star\\ M_2 & M_3
\end{bmatrix}.
\]
When $M_2=0$, we use $\blkdiag(M_1,M_3)$ to denote the resulting block-diagonal matrix. For $x\in\mathbb{R}^n$, we use $\norm*{x}$ to denote the Euclidean norm. We also denote the spectral norm of a matrix $A$ by $\norm{A}$.

In this paper, we will introduce a notion of descriptor embedding for nonlinear systems. Therefore, some preliminaries are required regarding the class of linear parameter-varying descriptor (LPVd) systems. Such systems take the form
\begin{equation}\label{eqn:polytopic_singularLPV}
	E \eta_{t+1} = A(\rho_t)\eta_t,
\end{equation}
where $\eta_t\in\mathbb{R}^{\nu}$ is the state vector, $E\in\mathbb{R}^{\nu\times \nu}$ is a singular matrix with $\rank(E)=r < \nu$, $\rho_t\in\mathcal{D}\subseteq\mathbb{R}^{p}$ is a time-varying parameter and $t\in\mathbb{N}$. The matrix $A(\rho_t)$ is assumed to belong to a polytope
\begin{align*}
	\mathcal{A} \coloneqq 
	\Big\lbrace A(\rho_t) ~\Big|~ 
	A(\rho_t) = \sum\nolimits_{i=1}^N \alpha_i(\rho_t) A_i,\quad A_i\in\mathbb{R}^{\nu\times \nu},\\
	\alpha_i(\rho_t) \geq 0,\quad\sum\nolimits_{i=1}^N \alpha_i(\rho_t) = 1
	\Big\rbrace.
\end{align*}
Throughout this work, we make the mild assumption that the functions $\alpha_i$ are continuous. The following is an important property of LPVd systems.
\begin{definition}\label{def:polyquad_admiss}
	System \eqref{eqn:polytopic_singularLPV} is \textit{poly-quadratically admissible} if it is regular and causal\footnote{Regularity guarantees existence and uniqueness of solutions for \eqref{eqn:polytopic_singularLPV}. Causality is in the usual sense that the system's state at time $t$ depends solely on inputs up to time $t$, see \cite{bara2011dilated,barbosa2018} and references therein for more details.} over the polytope $\mathcal{A}$, and there exists a parameter-dependent Lyapunov function (PDLF) $V(\eta)=~\eta^\top E^\top P(\rho) E \eta$ satisfying a strict decrease condition $V(\eta_{t+1}) - V(\eta_t) < 0,\,\forall\eta\neq0$, where $P(\rho_t) = \sum\nolimits_{i=1}^N \alpha_i(\rho_t) P_i,\, P_i\succ0$.\hfill$\square$
\end{definition}
~
Necessary and sufficient conditions for poly-quadratic admissibility are provided in the so-called singular-value decomposition (SVD) canonical form, which is an equivalent representation of system \eqref{eqn:polytopic_singularLPV}. In particular, let
\[
R E W = \begin{bmatrix}
	I_r &0_{\nu-r\times r} \\ 0_{r\times\nu-r} & 0_{\nu-r\times \nu-r}
\end{bmatrix},
\]
which can be obtained from the singular value decomposition of the matrix $E$, where $R= \blkdiag(\Sigma^{-1},I_{\nu-r})U^\top$ with $\Sigma$ being a diagonal matrix containing the non-zero singular values of $E$, whereas $U,W$ are unitary matrices containing the left/right singular vectors, respectively. Now partition $RA(\rho_t)W$ accordingly as
\[
\begin{aligned}
	R A(\rho_t) W &= \begin{bmatrix}
		A_{11}(\rho_t) & A_{12}(\rho_t)\\ A_{21}(\rho_t) & A_{22}(\rho_t)
	\end{bmatrix},
\end{aligned}
\]
and define a coordinate transformation $\eta_{t} = W \bar{\eta}_t$ with $\bar{\eta}_t = \begin{bmatrix}
	\bar{\eta}_{1,t}^\top & \bar{\eta}_{2,t}^\top
\end{bmatrix}^\top$ where $\bar{\eta}_{1,t}\in\mathbb{R}^{r}$ and $\bar{\eta}_{2,t}\in\mathbb{R}^{\nu-r}$. Then, by substituting $\eta_t = W\bar{\eta}_t$ in \eqref{eqn:polytopic_singularLPV} and pre-multiplying by $R$, one obtains the following equivalent SVD canonical form
\begin{equation}\label{eqn:polytopic_singularLPV_SVDform}
	\begin{bmatrix}
		I_r & 0\\0 & 0
	\end{bmatrix}\begin{bmatrix}
		\bar{\eta}_{1,t+1}\\ \bar{\eta}_{2,t+1}
	\end{bmatrix} = \begin{bmatrix}
		A_{11}(\rho_t) & A_{12}(\rho_t)\\ A_{21}(\rho_t) & A_{22}(\rho_t)
	\end{bmatrix}\begin{bmatrix}
		\bar{\eta}_{1,t}\\ \bar{\eta}_{2,t}
	\end{bmatrix}.
\end{equation}

Notice that in these coordinates, the system is split into two parts: a dynamic subsystem $\bar{\eta}_1$ and algebraic constraints $A_{21}(\rho_t)\bar{\eta}_{1,t} + A_{22}(\rho_t)\bar{\eta}_{2,t} = 0$. The following theorem provides necessary and sufficient conditions for poly-quadratic admissibility of system \eqref{eqn:polytopic_singularLPV}.
\begin{theorem}[\!\!{\cite{bara2011dilated,barbosa2018}}]\label{thm:SVDform}
	System \eqref{eqn:polytopic_singularLPV} (or equivalently \eqref{eqn:polytopic_singularLPV_SVDform}) is poly-quadratically admissible if and only if $A_{22}(\rho_t)$ is invertible for all $\rho_t\in\mathcal{D}\subseteq\mathbb{R}^{p}$, and 
	\begin{equation}
		\begin{aligned}
			\bar{\eta}_{1,t+1} &= \left( A_{11}(\rho_t) - A_{12}(\rho_t) A_{22}^{-1}(\rho_t) A_{21}(\rho_t) \right) \bar{\eta}_{1,t}\\
			&\eqqcolon \bar{A}(\rho_t)\bar{\eta}_{1,t} 
		\end{aligned}
		\label{eqn:bareta1}
	\end{equation}
	is \textit{poly-quadratically stable}.\hfill$\square$
\end{theorem}

Poly-quadratic stability is a notion of stability originally introduced in \cite[Def. 2]{daafouz2001parameter} for polytopic \emph{linear} systems. A system of the form \eqref{eqn:bareta1} is said to be poly-quadratically stable if there exists a PDLF $\bar V(\bar{\eta}_1) = \bar{\eta}_1^\top \bar{P}(\rho) \bar{\eta}_1$, for some $\bar{P}(\rho_t)=\sum_{i=1}^{N}\alpha_i(\rho_t)\bar{P}_i$, $\bar{P}_i\succ0\in\mathbb{R}^{r\times r}$, satisfying a strict decrease condition\footnote{As a result, poly-quadratic stability implies that the origin of \eqref{eqn:bareta1} is uniformly asymptotically stable for all $\rho\in\mathcal{D}\subseteq\mathbb{R}^{p}$, compare~\cite{daafouz2001parameter}.} $\bar{V}(\bar\eta_{1,t+1})-\bar{V}(\bar\eta_{1,t})<0$ for all $\bar\eta_1\neq0$, i.e., that the following inequality is satisfied for all $\rho\in\mathcal{D}\subseteq\mathbb{R}^{p}$
\begin{equation}
	\bar{A}(\rho_{t})^\top \bar{P}(\rho_{t+1}) \bar{A}(\rho_{t}) - \bar{P}(\rho_t) \prec 0.
	\label{eqn:PDLFinequalitybareta1}
\end{equation}

Although Theorem~\ref{thm:SVDform} provides necessary and sufficient conditions for admissibility of \eqref{eqn:polytopic_singularLPV}, \textit{testing} for poly-quadratic stability of \eqref{eqn:bareta1} is not always possible using the techniques developed in \cite{daafouz2001parameter}. This is because \eqref{eqn:bareta1} does not necessarily have a polytopic structure (due to the term $A_{12}(\rho_t) A_{22}^{-1}(\rho_t) A_{21}(\rho_t)$, compare \cite[Remark 2]{de2008robust}). As an alternative, \cite{bara2011dilated} provides sufficient LMI conditions for poly-quadratic admissibility. This is summarized in the following theorem.
\begin{theorem}[\!\!{\cite[Thm. 4.3(ii)]{bara2011dilated}}]
	\label{thm:LPV_polyquadadmissibility}
	System \eqref{eqn:polytopic_singularLPV} (or equivalently \eqref{eqn:polytopic_singularLPV_SVDform}) is poly-quadratically admissible if there exist symmetric positive definite matrices $Q_i\succ0\in\mathbb{R}^{r\times r}$ for $i\in\mathcal{I}\coloneqq\{1,\ldots,N\}$, matrices $M,F\in\mathbb{R}^{r\times r}$, $C_1,C_3\in\mathbb{R}^{\nu-r\times r}$ and $C_2\in\mathbb{R}^{\nu-r\times \nu-r}$ such that
	\begin{equation}
		\resizebox{\columnwidth}{!}{$\begin{bmatrix}
				\sym\left( RA_iW\begin{bmatrix}
					F & 0\\ C_3 & C_2
				\end{bmatrix} \right) + \begin{bmatrix}
					Q_j & 0 \\ 0 & 0
				\end{bmatrix} & \star \\
				\left( RA_iW\begin{bmatrix}
					M\\C_1
				\end{bmatrix}+\begin{bmatrix}
					F^\top \\ 0
				\end{bmatrix} \right)^\top & \sym(M) - Q_i
			\end{bmatrix}\succ0,$}
		\label{eqn:polyquadraticadmissibilityLMIconditions}
	\end{equation}
	for all $(i,j)\in\mathcal{I}\times\mathcal{I}$. In this case, the PDLF of the corresponding subsystem \eqref{eqn:bareta1} takes the form $\bar V(\bar{\eta}_1) = \bar{\eta}_1^\top \bar{P}(\rho) \bar{\eta}_1$ with $\bar P(\rho)=\sum_{i=1}^N \alpha_i(\rho)Q_i^{-1}$.\hfill$\square$
\end{theorem}

\begin{remark}\label{rem:PDLFLPVd}
	In \cite{bara2011dilated}, the structure of the PDLF for the subsystem \eqref{eqn:bareta1} was not explicitly stated. However, it is easy to verify that it takes the form written in Theorem~\ref{thm:LPV_polyquadadmissibility}. For completeness, we include this in Appendix~\ref{app:onPDLF}. \hfill$\square$
\end{remark}

Now consider an LPVd system with inputs
\begin{equation}
	E\eta_{t+1} = A(\rho_t)\eta_t + B u_t,
	\label{eqn:polytopic_singularLPV_withinputs}
\end{equation}
where $u_t\in\mathbb{R}^m$ is the input at time $t$ and $B\in\mathbb{R}^{\nu\times m}$. For such systems, one is not only interested in stabilization but in \textit{admissibilization}. This is because properties of regularity and causality are not feedback invariant and need to be ensured, along with stability, for the closed-loop system. The problem of \textit{poly-quadratic admissibilization} was also addressed in \cite{bara2011dilated}. Specifically, the vertices $\kappa_i$ of a polytopic state-feedback controller $u=\kappa(\rho)\eta$, with
\[
\kappa(\rho) = \sum\nolimits_{i=1}^N \alpha_i(\rho) \kappa_i,\quad \kappa_i\in\mathbb{R}^{m\times \nu},
\]
are designed such that the closed-loop system is poly-quadratically admissible.

In the following section, we present a notion of \textit{descriptor embedding of nonlinear systems}, and later exploit Theorem~\ref{thm:LPV_polyquadadmissibility} to design stabilizing controllers for nonlinear systems.
\section{Descriptor embedding for nonlinear systems}\label{sec:embedding}
We consider discrete-time nonlinear systems of the form\footnote{For clarity of presentation, we first present the results for systems with constant input vector fields. Later in Section~\ref{sec:inaff} we discuss a more general class of nonlinear systems.}
\begin{equation}\label{eqn:NLsys}
    x_{t+1} = f(x_t)+Bu_t
\end{equation}
where $x_t\in\mathbb{R}^n$, $u_t\in\mathbb{R}^m$ are the state and input vectors at time $t\in\mathbb{N}$, respectively. Both $f:\mathbb{R}^n\to\mathbb{R}^n$ and the matrix $B\in\mathbb{R}^{n\times m}$ are unknown, but $f$ satisfies the following assumption.

\begin{assumption}\label{asmp:f}
    The function $f$ is continuous, zero at zero, and can be written as $f(x)=AZ(x)$ where $A\in\mathbb{R}^{n\times S}$ is an unknown matrix and $Z:\mathbb{R}^n\to\mathbb{R}^{S}$ is a known, continuous, vector-valued function that takes the form $Z(x) = \begin{bmatrix}
        x^\top & \underline{Z}(x)^\top
    \end{bmatrix}^\top$, with $Z(0)=0$ and $\underline{Z}:\mathbb{R}^n\to\mathbb{R}^{S-n}$.\hfill$\square$
\end{assumption}
\begin{remark}\label{rem:Z}
    The user-defined basis functions $Z$ may be chosen, e.g., based on physical knowledge about the system. Requiring that $Z(0)=0$ is necessary for the subsequent results on stabilizing the origin of the nonlinear system \eqref{eqn:NLsys}. This is, however, not restrictive as basis functions that do not satisfy this assumption can be suitably shifted by a constant. Later in Section~\ref{sec:inexact}, we relax Assumption~\ref{asmp:f} and account for basis functions approximation error.\hfill$\square$
\end{remark}
~
Assumption~\ref{asmp:f} allows us to write the system \eqref{eqn:NLsys} in the following form
\begin{align}
	x_{t+1} &= AZ(x_t) + Bu_t.\label{eqn:NLsys_Z}
\end{align}
We seek a polytopic state-feedback controller $u=K(x)Z(x)$ such that the origin of the closed-loop nonlinear system
\begin{equation}
	x_{t+1} = (A+BK(x_t))Z(x_t) \label{eqn:CL_NLsys}
\end{equation}
is asymptotically stable. By polytopic state-feedback control, we mean that $K(x)$ takes the form
\begin{equation}
	K(x) = \sum\nolimits_{i=1}^N \alpha_i(x)K_i, \quad K_i\in\mathbb{R}^{m\times S}, \label{eqn:CLgain}
\end{equation}
where $\alpha_i$ are known\footnote{These functions will be specified later, see Assumption~\ref{asmp:L} below.}, continuous, non-negative (i.e., $\alpha_i(x)\geq0$) functions that partition to unity (i.e., $\sum\nolimits_{i=1}^N\alpha_i(x)=1$). Since $A,B$ are unknown, we perform an experiment on the open-loop system \eqref{eqn:NLsys} and collect input-state data $\{u_k^d,x_k^d\}_{k=0}^{T}$ that are arranged in the following matrices
\begin{equation}\label{eqn:data_mats}
	\begin{aligned}
		U_0 &= \begin{bmatrix}
			u_0^d & u_1^d & \cdots & u_{T-1}^d
		\end{bmatrix}\in\mathbb{R}^{m\times T},\\
		Z_0 &= \begin{bmatrix}
			Z(x_0^d) & Z(x_1^d) & \cdots & Z(x_{T-1}^d)
		\end{bmatrix}\in\mathbb{R}^{S\times T},\\
		X_1 &= \begin{bmatrix}
			x_1^d & x_2^d & \cdots & x_{T}^d
		\end{bmatrix}\in\mathbb{R}^{n\times T},  
	\end{aligned}
\end{equation}
and satisfy the following equation (compare \eqref{eqn:NLsys_Z})
\begin{equation}\label{eqn:data_identity}
	X_1 = AZ_0 + BU_0.
\end{equation}
The following lemma provides a data-dependent representation of the closed-loop system. The results of this lemma are analogous to, e.g., \cite{Guo2022polynomial,DePersis2023}, with the main difference being the use of a polytopic state-feedback gain $K(x)$.
~
\begin{lemma}\label{lemma:DBrep}
	Consider matrices $K_i\in\mathbb{R}^{m\times S}, G_i\in\mathbb{R}^{T\times S}$, for $i\in\mathcal{I}\coloneqq\{1,\ldots,N\}$, satisfying
	\begin{equation}
		\begin{bmatrix}
			U_0\\ Z_0
		\end{bmatrix}G_i = \begin{bmatrix}
			K_i\\ I_{S}
		\end{bmatrix}.
		\label{eqn:KGrelation}
	\end{equation}
	Then, \eqref{eqn:CL_NLsys} can be equivalently written as
	\begin{equation}
		x_{t+1} = X_1G(x_t)Z(x_t) \label{eqn:CL_NLsys_data}
	\end{equation}
	where $G(x) = \sum\nolimits_{i=1}^{N}\alpha_i(x)G_i$ and $\alpha_i$ defined below \eqref{eqn:CLgain}.\hfill$\square$
\end{lemma}
\begin{proof}
	Let \eqref{eqn:KGrelation} hold. Using the definition of $K(x)$ in \eqref{eqn:CLgain}, we have
	\begin{align}
		K(x) = \sum\nolimits_{i=1}^{N}\alpha_i(x)K_i &\stackrel{\eqref{eqn:KGrelation}}{=} \sum\nolimits_{i=1}^{N}\alpha_i(x)U_0G_i\label{eqn:Kx_U0Gx}\\
		&= U_0 \sum\nolimits_{i=1}^{N}\alpha_i(x)G_i = U_0G(x).\notag
	\end{align}
	Furthermore, since $Z_0G_i = I_{S}$ holds for all $i\in\mathcal{I}$, one can multiply both sides of $Z_0G_i = I_{S}$ by $\alpha_i(x)$ which gives $\alpha_i(x) Z_0G_i = \alpha_i(x)I_{S}$. Summing these $N$ equations gives
	\begin{equation}
		\begin{aligned}
			\sum\nolimits_{i=1}^{N}\alpha_i(x) Z_0 G_i &= \sum\nolimits_{i=1}^{N}\alpha_i(x) I_{S}\\
			\Rightarrow Z_0 \sum\nolimits_{i=1}^{N}\alpha_i(x)G_i &= I_{S}\\
			\Rightarrow Z_0 G(x) &= I_{S}
		\end{aligned}\label{eqn:Z0Gx_I}
	\end{equation}
	where we have made use of $\sum_{i=1}^{N}\alpha_i(x)=1$. Equations~\eqref{eqn:Kx_U0Gx} and \eqref{eqn:Z0Gx_I} together give
	\begin{equation}\label{eqn:KGxrelation}
		\begin{bmatrix}
			U_0 \\ Z_0
		\end{bmatrix}G(x) = \begin{bmatrix}
			K(x)\\ I_{S}
		\end{bmatrix}.
	\end{equation}
	Finally, one can write $A+BK(x)$ in \eqref{eqn:CL_NLsys} as
	\[
	\begin{aligned}
		A+BK(x) &= \begin{bmatrix}
			B & A
		\end{bmatrix}\begin{bmatrix}
			K(x)\\ I_{S}
		\end{bmatrix} \\ &\stackrel{\eqref{eqn:KGxrelation}}{=} \begin{bmatrix}
			B & A
		\end{bmatrix}\begin{bmatrix}
			U_0 \\ Z_0
		\end{bmatrix}G(x) \stackrel{\eqref{eqn:data_identity}}{=} X_1G(x),
	\end{aligned}
	\]
	which completes the proof.
\end{proof}%
\begin{remark}
	A necessary condition for \eqref{eqn:KGrelation} to hold is that $Z_0$ has full row rank, which can be viewed as a sufficient richness condition that can be easily verified from recorded data. This is weaker than persistence of excitation conditions that are typically employed in majority of recent works on direct data-driven control, compare \cite{markovsky2021behavioral,martin2023guarantees,depersis2023annrev} and references therein. Under certain conditions, it can be enforced by design of input, see \cite{Alsalti_pe}. Later, we will see that $Z_0$ having full row rank is necessary for the feasibility of certain data-dependent LMIs whose solution yields a stabilizing controller.\hfill$\square$
\end{remark}

Lemma~\ref{lemma:DBrep} provides a data-based representation of the closed-loop system. This now transforms the problem of stabilization to one of finding $G(x)$ such that the origin of \eqref{eqn:CL_NLsys_data} is asymptotically stable. To do this, we first make the following assumption.
\begin{assumption}\label{asmp:L}
	There exists a function $L:\mathbb{R}^n\to\mathbb{R}^{S-n\times S}$ such that, for all $x\in\mathcal{X}\subseteq\mathbb{R}^n$ with $0\in\mathrm{int}(\mathcal{X})$,
	\begin{adjustwidth}{1.5em}{0pt}
		\begin{itemize}
			\item[(A2.1)] $L(x)Z(x)=0$,
			\item[(A2.2)] $L(x) = \begin{bmatrix}L_1(x) & L_2(x)\end{bmatrix}$ with $L_2(x)\in\mathbb{R}^{S-n\times S-n}$ invertible, and
			\item[(A2.3)] $L(x) = \sum\nolimits_{i=1}^N \alpha_i(x) \mathcal{L}_i$, where $\mathcal{L}_i\coloneqq L(v_i)$ for some $v_i\in\mathcal{X}$ and $\alpha_i$ are continuous user-defined functions satisfying $\alpha_i(x)\geq0$ and $\sum_{i=1}^{N}\alpha_i(x)=1$.\hfill$\square$
		\end{itemize}
	\end{adjustwidth}	
\end{assumption}
Later in this section, we discuss Assumption~\ref{asmp:L} in more detail and point out classes of $Z$ for which a function $L$ satisfying it exists (see Theorem~\ref{thm:asmpLmonomoialZ}). Intuitively, this assumption is a constraint on the system that must be satisfied at all times. Indeed, it holds by (A2.1) and (A2.2) that $\underline{Z}(x) = -L_2^{-1}(x)L_1(x)x$ and the closed-loop system can be written as
\begin{equation}
	\begin{aligned}
		x_{t+1} &= X_1G(x_t)Z(x_t)\\
		&= \left( X_1\overline{G}(x_t) - X_1\underline{G}(x_t)L_2^{-1}(x_t)L_1(x_t) \right)x_t,
	\end{aligned}
	\label{eqn:CL_NLsys_LPV}
\end{equation}
where $G(x) = \begin{bmatrix} \overline{G}(x) & \underline{G}(x)\end{bmatrix}$ is partitioned appropriately. Notice that this structure is reminiscent of that in~\eqref{eqn:bareta1}, with $\bar\eta_1 = \rho = x \in \mathcal{X}\subseteq\mathbb{R}^n$. The following theorem shows that poly-quadratic admissibility of a certain polytopic LPVd system implies asymptotic stability of the origin of \eqref{eqn:CL_NLsys_LPV}. We refer to this result as the \emph{descriptor embedding of nonlinear systems} and is the main result of this section.
\begin{theorem}\label{thm:embedding}
	Let Assumptions~\ref{asmp:f}-\ref{asmp:L} and condition \eqref{eqn:KGrelation} hold, and consider a polytopic LPVd system
	\begin{equation}
		\begin{bmatrix}
			I_n & 0 \\ 0 & 0
		\end{bmatrix}\eta_{t+1} = \begin{bmatrix}
			X_1G(\rho_t)\\ L(\rho_t)
		\end{bmatrix}\eta_t,
		\label{eqn:embedded_CLsys}
	\end{equation}
	for which $\rho_t \in \mathcal{X} \subseteq \mathbb{R}^n$ for all $t \in \mathbb{N}$, and $G$ as in Lemma~\ref{lemma:DBrep}. If \eqref{eqn:embedded_CLsys} is poly-quadratically admissible, then the origin of \eqref{eqn:CL_NLsys_LPV} is locally asymptotically stable. Furthermore, if $\mathcal{X}=\mathbb{R}^n$, then the origin is globally asymptotically stable.
\end{theorem}
\begin{proof}
	Notice that \eqref{eqn:embedded_CLsys} is in SVD canonical form (compare \eqref{eqn:polytopic_singularLPV_SVDform}). In particular,
	\[
	\begin{bmatrix}
		I_n & 0 \\ 0 & 0
	\end{bmatrix}\begin{bmatrix}
		\bar\eta_{1,t+1} \\ \bar\eta_{2,t+1}
	\end{bmatrix} = \begin{bmatrix}
		X_1 \overline{G}(\rho_t) & X_1\underline{G}(\rho_t) \\ L_1(\rho_t) & L_2(\rho_t)
	\end{bmatrix}\begin{bmatrix}
		\bar\eta_{1,t} \\ \bar\eta_{2,t}
	\end{bmatrix}
	\]
	where $\bar\eta_1$ denotes the first $n$ elements of $\eta$ while $\bar\eta_2$ denotes the remaining ones. Since $\rho_t \in \mathcal{X} \subseteq \mathbb{R}^n$ for all $t \in \mathbb{N}$, it holds by\footnote{Notice that (A2.2) is automatically satisfied if \eqref{eqn:embedded_CLsys} is poly-quadratically admissible, compare Theorem~\ref{thm:SVDform}. Nevertheless, we keep it as an assumption to make it clear to the user that $L$ needs to be designed such that (A2.2) holds, as otherwise the data-dependent LMIs presented later are not feasible.} (A2.2) that $L_2(\rho_t)$ is invertible and hence, by Theorem~\ref{thm:SVDform}, poly-quadratic admissibility of \eqref{eqn:embedded_CLsys} is equivalent to the following dynamic subsystem being poly-quadratically stable
	\begin{equation}
		\bar\eta_{1,t+1} = \left( X_1\overline{G}(\rho_t) - X_1\underline{G}(\rho_t)L_2^{-1}(\rho_t)L_1(\rho_t) \right)\bar\eta_{1,t},
		\label{eqn:embedded_subsystem}
	\end{equation}
	with a corresponding Lyapunov function of the form\footnote{If poly-quadratic admissibility of \eqref{eqn:embedded_CLsys} is verified using, e.g., Theorem~\ref{thm:LPV_polyquadadmissibility}, then one can take $\bar{P}_i=Q_i^{-1}$, compare Remark~\ref{rem:PDLFLPVd}.} $ V(\bar\eta_1) = \bar\eta_1^\top \bar{P}(\rho) \bar\eta_1$, $\bar{P}(\rho)=\sum_{i=1}^N \alpha_i(\rho)\bar{P}_i$, satisfying a strict decrease condition $V(\bar\eta_{1,t+1})- V(\bar\eta_{1,t})<0$ for all $\bar\eta_1\neq0$. Using standard Lyapunov arguments \cite{haddad2008nonlinear} and the fact that $\mathcal{X}$ has a non-empty interior (see Assumption~\ref{asmp:L}), there exists a set $\Omega\subseteq\mathcal{X}$ with $0\in\mathrm{int}(\Omega)$ such that $\bar\eta_{1,t}\in\Omega\implies\bar\eta_{1,t+1}\in\Omega$, i.e., $\Omega$ is positively invariant for system~\eqref{eqn:embedded_subsystem}. 
	
	Now consider the following system
	\begin{equation}
		z_{t+1} = \left( X_1\overline{G}(r_t) - X_1\underline{G}(r_t)L_2^{-1}(r_t)L_1(r_t) \right)z_t
		\label{eqn:zdynamics_rt}
	\end{equation}
	with some parameter
	\begin{equation}
		r_t \coloneqq \begin{cases}
			z_t, \quad & z_t\in\Omega,\\
			\mathrm{Proj}_{\Omega}(z_t), \quad & \mathrm{otherwise},
		\end{cases}
		\tag{22}
	\end{equation}
	where $\mathrm{Proj}_{\Omega}(\mu)\coloneqq \arg\min\limits_{y\in\Omega} \norm{y-\mu}_2$ is the projection of $\mu$ on the set $\Omega$ (which exists since $\Omega$ is closed). Notice that by construction of $r_t$ it holds that $r_t\in\Omega\subseteq\mathcal{X}$ for all time (regardless of how $z_t$ evolves). Moreover, system \eqref{eqn:zdynamics_rt} takes the same form as \eqref{eqn:embedded_subsystem} but with a different parameter (which, by construction, remains in $\Omega\subseteq\mathcal{X}$ for all time). Therefore, from the above arguments it follows that $\Omega$ is also a forward invariant set for \eqref{eqn:zdynamics_rt}. Now consider $z_0\in\Omega$, then $r_0=z_0\in\Omega$. By forward invariance of $\Omega$ for \eqref{eqn:zdynamics_rt}, we have $z_{1}\in\Omega$ and, hence, $r_1=z_1\in\Omega$. By induction, we find that $z_0\in\Omega\implies r_t\equiv z_t \in\Omega\subseteq\mathcal{X}$ for all $t\in\mathbb{N}$. Said differently, if $z_0\in\Omega$, then $z$ evolves according to
	\[
	z_{t+1} = \left( X_1\overline{G}(z_t) - X_1\underline{G}(z_t)L_2^{-1}(z_t)L_1(z_t) \right)z_t
	\]
	which is exactly \eqref{eqn:CL_NLsys_LPV}. Moreover, $V(z)=z^\top \bar{P}(z) z$, with $\bar{P}(z)$ as above, is a Lyapunov function for the system satisfying $V(z_{t+1})-V(z_t)<0$ for all $z\in\Omega\backslash\{0\}$. Asymptotic stability then directly follows from standard Lyapunov arguments \cite{haddad2008nonlinear}. Finally, since $V$ is radially unbounded, it follows that the results are global if $\mathcal{X}=\mathbb{R}^n$ (in which case $\Omega=\mathcal{X}$).
\end{proof}

Theorem~\ref{thm:embedding} forms the basis of the subsequent developments and deserves particular attention. Most importantly, it allows us to infer stability of the nonlinear system by enforcing poly-quadratic admissibility of some higher-dimensional descriptor system. In particular, if there exist $G_i$, for all $i\in\mathcal{I}$, such that~\eqref{eqn:embedded_CLsys} is poly-quadratically admissible, then a controller of the form~\eqref{eqn:CLgain} asymptotically stabilizes the origin of the closed-loop nonlinear system~\eqref{eqn:CL_NLsys_LPV} (which is equivalent to \eqref{eqn:CL_NLsys_data} or \eqref{eqn:CL_NLsys} on $\mathcal{X}$). In the following section, we exploit Theorems~\ref{thm:LPV_polyquadadmissibility} and \ref{thm:embedding} to derive data-dependent LMIs which, if feasible, return matrices $G_i$ such that \eqref{eqn:embedded_CLsys} is poly-quadratically admissible. The stabilizing controller for the nonlinear system can then be obtained as $K(x)=U_0G(x)$, see Lemma~\ref{lemma:DBrep}. 

\begin{remark}
	Theorem~\ref{thm:embedding} focused on using tools from descriptor systems theory to infer stability of nonlinear systems satisfying Assumptions~\ref{asmp:f} and \ref{asmp:L}. We expect that the insights gained from this novel embedding notion may lead to new nonlinear system analysis and control design methodologies (both in model-based and data-based contexts) beyond what is presented in this paper. Future research will focus on leveraging other tools from descriptor systems theory to further analyze and control nonlinear systems. \hfill$\square$
\end{remark}

We conclude this section with the following discussion concerning Assumption~\ref{asmp:L}. Although it may potentially be restrictive, we prove in the following theorem that such a function $L$ exists when $Z$ is composed of all monomials of the state $x$ up to some finite degree and $\mathcal{X}$ is a polytope (i.e., convex and bounded polyhedron).
~
\begin{theorem}\label{thm:asmpLmonomoialZ}
	Consider a function $Z:\mathbb{R}^n\to\mathbb{R}^S$ of the form $Z(x) = \begin{bmatrix}x^\top & \underline{Z}(x)^\top\end{bmatrix}^\top$ where $\underline{Z}:\mathbb{R}^n\to\mathbb{R}^{S-n}$ consists of all monomials up to some finite degree $t\geq 2$. Then, there exists a function $L:\mathcal{X}\to\mathbb{R}^{S-n\times S}$ which satisfies Assumption~\ref{asmp:L} with $\mathcal{X}\subset\mathbb{R}^n$ being a polytope.\hfill$\square$
\end{theorem}
\begin{proof}
	We will show the claim in two steps. First, we will construct $L$ by making use of the fact that monomials can be factorized as products of its lower degree terms. Then, we show that the constructed $L$ satisfies Assumption~\ref{asmp:L}. 
	
	To start, let the elements of $\underline{Z}$ be denoted by $\underline{z}_j(x)$, for $j\in\{1,\ldots,S-n\}$, and ordered in ascending manner according to their degrees, i.e., $\deg(\underline{z}_1(x))\leq\deg(\underline{z}_2(x))\leq\ldots\leq\deg(\underline{z}_{S-n}(x))$. Consider the first monomial term $\underline{z}_1(x)$ in $\underline{Z}(x)$ which is of degree~2. There exists a \textit{linear} vector-valued function $\ell_1:\mathbb{R}^n\to\mathbb{R}^{n}$ such that\footnote{For clarity, the following procedure is illustrated in Appendix~\ref{app:exmpL} on an example of all monomials up to degree three for $x\in\mathbb{R}^2$.}
	\[
	\left[\begin{array}{c|cccc}
		\ell_1^\top(x) & -1 & 0 & \cdots & 0
	\end{array}\right]\begin{bmatrix}
	x\\ \hline \underline{z}_1(x) \\ \underline{z}_{2}(x) \\ \vdots \\ \underline{z}_{S-n}(x)
\end{bmatrix} = 0.
	\]
	Similarly, for the second entry $\underline{z}_2(x)$ of $\underline{Z}(x)$, there exists a linear function $\ell_2:\mathbb{R}^n\to\mathbb{R}^{n+1}$ such that
	\[
	\left[\begin{array}{c|cccc}
		\ell_2^\top(x) & -1 & 0 & \cdots & 0
	\end{array}\right]\begin{bmatrix}
		x\\ \underline{z}_1(x) \\\hline \underline{z}_{2}(x) \\ \vdots \\ \underline{z}_{S-n}(x)
	\end{bmatrix} = 0.
	\]
	Proceeding in the same manner, we have that for the $j$th entry $\underline{z}_j(x)$ of $\underline{Z}(x)$, for any $j\in\{1,\ldots,S-n\}$, there exists a linear function $\ell_j:\mathbb{R}^n\to\mathbb{R}^{n+j-1}$ (note that each $\ell_j$ maps to spaces of different dimensions) such that 
	\[
	\big[\ell_j^\top(x) ~\big|~ -1 ~\,\, \underbrace{\begin{matrix}
					0 & \cdots & 0
		\end{matrix}}_{S-n-j}\big]\begin{bmatrix}
		x\\ \vdots \\ \hline \underline{z}_j(x) \\ \underline{z}_{j+1}(x)\\ \vdots\\ \underline{z}_{S-n}(x)
	\end{bmatrix} = 0.
	\]
	Concatenating all the above results in
	\[
		\begin{tikzpicture}
			\draw[draw=black] (-0.78,-1.4) -- (-0.78,-1.25);
			\draw[draw=black] (-0.78,-0.725) -- (-0.78,-0.13);
			\draw[draw=black] (-0.78,0.3) -- (-0.78,0.75);
			\draw[draw=black] (1.0,-1.0175) -- (2,-1.0175);
			\node at (0,0) {
				\begingroup
				\setlength{\arraycolsep}{8pt}
				$\begin{aligned}
					L(x)&\coloneqq \left[\begin{array}{c|c}
						L_1(x) & L_2(x)
					\end{array}\right]\\
					&=\left[\,\begin{array}{ccccc}
						\ell_1^\top(x) & -1 & 0 & \cdots & 0 \\
						\multicolumn{2}{c}{\rule[.5ex]{1em}{0.4pt}\,\, \ell_{2}^\top(x) \,\, \rule[.5ex]{1em}{0.4pt}} & -1 & \cdots & 0\\
						& & & \ddots & \vdots\\
						\multicolumn{3}{c}{\rule[.5ex]{2em}{0.4pt}\,\,\ell_{S-n}^\top(x)\,\,\rule[.5ex]{2em}{0.4pt}} & &-1
					\end{array}\right].
				\end{aligned}$
				\endgroup
				};
		\end{tikzpicture}
	\]
	Notice that, by construction, $L(x)Z(x)=0$ and hence satisfies Assumption (A2.1). Moreover, note that $L_2(x)$ has a lower-triangular structure with the diagonal entries being $-1$ and, hence, $L_2(x)$ is invertible for all $x$ thus satisfying (A2.2). Finally, notice that the resulting function $L$ is affine in $x$ (since each $\ell_j$ is of degree 1). Since, moreover, $\mathcal{X}\subset\mathbb{R}^n$ is a polytope, it follows that (A2.3) holds since the image of a polytope (i.e., convex and bounded polyhedron) under an affine map of is also a polytope, see \cite[p. 384]{toth2017handbook}.
\end{proof}

The fact that Assumption~\ref{asmp:L} is satisfied when choosing $Z$ as monomials implies that our subsequent analysis is directly applicable to the class of polynomial systems, albeit only locally, since Theorem~\ref{thm:asmpLmonomoialZ} only guarantees that Assumption~\ref{asmp:L} holds on (arbitrarily large but bounded) polytopes $\mathcal{X}\subset\mathbb{R}^n$. Additionally, choosing $Z$ as monomials is advantageous due to their ability to tightly approximate smooth vector fields $f$ on compact regions \cite{Stone1948}. Later in Section~\ref{sec:inexact}, we account for basis functions approximation errors in the context of the proposed framework. 

For the case of $L$ as in Theorem~\ref{thm:asmpLmonomoialZ}, one can always find continuous functions $\alpha_i$ satisfying (A2.3) using the generalized barycentric coordinates \cite{floater2015generalized} with respect to $\mathcal{X}$, see also Appendix~\ref{app:exmpL}. Knowledge of $\alpha_i$ is important as they are used in the polytopic state-feedback controller, compare \eqref{eqn:CLgain}. For general choices of $L$, the functions $\alpha_i$ are typically determined on a case-by-case basis, as no systematic method is readily available.

Finally, we note that Assumption~\ref{asmp:L} may hold globally for some other (non-polynomial) choices of $Z$. This is illustrated in Example~\ref{ex:global_inv_pend} below, where, moreover, the resulting stabilizing controller is a global one.
\section{Data-based stabilization of nonlinear systems}
\label{sec:nominal}
\subsection{Main result}
In this section, we aim to determine matrices $G_i$ such that the system \eqref{eqn:embedded_CLsys} is poly-quadratically admissible. In particular, we make use of Theorems~\ref{thm:LPV_polyquadadmissibility} and \ref{thm:embedding} to propose the following result.

\begin{theorem}\label{thm:Kx}
    Consider a nonlinear system as in \eqref{eqn:NLsys} and let Assumptions~\ref{asmp:f}-\ref{asmp:L} hold. Fix some scalar $\delta\in\mathbb{R}$ and consider \eqref{eqn:P1} (shown at the top of the next page) in the decision variables $R_i\in\mathbb{R}^{T\times S},\,Q_i=Q_i^\top\succ0\in\mathbb{R}^{n\times n}$, for all $i\in\mathcal{I}=\{1,\ldots,N\}$, $M\in\mathbb{R}^{n\times n},\,C_1\in\mathbb{R}^{S-n \times n}$ and an invertible\footnote{We do not explicitly enforce invertibility of $C_2$ in \eqref{eqn:P1}. This is because, similar to \cite{bara2011dilated,masubuchi1997h}, if \eqref{eqn:P1} is feasible with a singular $C_2$, then the solutions $Q_i,M,C_1,C_2,R_i$ can be used to obtain a scalar $\epsilon\in(0,1)$ such that \eqref{eqn:P1} holds with $C_2$ replaced by $C_2+\epsilon I_{S-n}$.} $C_2\in\mathbb{R}^{S-n \times S-n}$.
    \begin{figure*}
        \centering
        \begin{subequations}
        \begin{align}
            &\begin{bmatrix}
                \sym\left(\begin{bmatrix}
                    X_1\\ \mathcal{L}_i Z_0
                \end{bmatrix}R_i\begin{bmatrix}
                    \delta I_{n} & 0\\0 & I_{S-n}
                \end{bmatrix}\right) + \begin{bmatrix}
                    Q_j & 0\\0 & 0
                \end{bmatrix} & \star \\ 
                \left(\begin{bmatrix}
                    X_1\\ \mathcal{L}_i Z_0
                \end{bmatrix}R_i\begin{bmatrix}
                    I_{n}\\ 0
                \end{bmatrix}+\begin{bmatrix}
                    \delta M^\top\\ 0
                \end{bmatrix}\right)^\top & \sym(M)-Q_i
            \end{bmatrix}\succ0,\quad\forall (i,j)\in\mathcal{I}\times\mathcal{I},\label{eqn:P1a}\\
            &Z_0 R_i = \begin{bmatrix}
                M & 0 \\ C_1 & C_2
            \end{bmatrix}.\label{eqn:P1b}
        \end{align}\label{eqn:P1}
    \end{subequations}
        \hrule
    \end{figure*}
    If \eqref{eqn:P1} is feasible, then the origin of \eqref{eqn:CL_NLsys}, with $K(x)=\sum_{i=1}^{N}\alpha_i(x)K_i$ and
    \begin{equation}
        K_i = U_0 R_i \begin{bmatrix}
            M & 0 \\ C_1 & C_2
        \end{bmatrix}^{-1},\label{eqn:P1_gain}
    \end{equation}
    is locally asymptotically stable. A Lyapunov function for the closed-loop system is given by $V(x_t)=x_t^\top P(x_t)x_t$ with $P(x_t) = \sum_{i=1}^{N}\alpha_i(x_t)Q_i^{-1}$. Furthermore, if $\mathcal{X}=\mathbb{R}^n$, then the origin is globally asymptotically stable.\hfill$\square$
\end{theorem}
\begin{proof}
	Let \eqref{eqn:P1} hold. From the bottom right block of \eqref{eqn:P1a}, it holds that $\sym(M)-Q_i\succ0$ for all $i\in\mathcal{I}$, and since $Q_i\succ0$ by assumption, it follows that $M$ is invertible. Furthermore, since $C_2$ is invertible by assumption, then one can define
	\begin{equation}
		G_i \coloneqq R_i\begin{bmatrix}
			M & 0\\ C_1 & C_2
		\end{bmatrix}^{-1}.\label{eqn:defGi}
	\end{equation}
	Pre-multiplying both sides from the left by $Z_0$ results in 
	\begin{equation}
		Z_0G_i = Z_0R_i\begin{bmatrix}
			M & 0\\ C_1 & C_2
		\end{bmatrix}^{-1}\stackrel{\eqref{eqn:P1b}}{=} I_{S},
	\label{eqn:enforced_Z0Gi}
	\end{equation}
	which, together with \eqref{eqn:P1_gain}, enforces \eqref{eqn:KGrelation} and, hence, the data-based representation of the closed-loop system \eqref{eqn:CL_NLsys_data} holds (see Lemma~\ref{lemma:DBrep}). 
	
	The proof proceeds by showing that \eqref{eqn:P1a} implies that the following LMIs hold
	\begin{equation}
		\resizebox{\columnwidth}{!}{$
			{\begin{bmatrix}
					\sym\left( \begin{bmatrix}
						X_1G_i\\ \mathcal{L}_i
					\end{bmatrix}\begin{bmatrix}
						F & 0\\ C_3 & C_2
					\end{bmatrix} \right) + \begin{bmatrix}
						Q_j & 0 \\ 0 & 0
					\end{bmatrix} & \star \\
					\left( \begin{bmatrix}
						X_1G_i\\ \mathcal{L}_i
					\end{bmatrix}\begin{bmatrix}
						M\\C_1
					\end{bmatrix}+\begin{bmatrix}
						F^\top \\ 0
					\end{bmatrix} \right)^\top & \sym(M) - Q_i
				\end{bmatrix}\succ0,}
			$}
		\label{eqn:polyquadadmissibility_proof}
	\end{equation}
	for some $F,C_3$ and for all $(i,j)\in\mathcal{I}\times\mathcal{I}$, which by Theorem~\ref{thm:LPV_polyquadadmissibility} implies poly-quadratic admissibility of some polytopic LPVd system of the form \eqref{eqn:embedded_CLsys} (which is already in SVD form) for any $\rho_t\in\mathcal{X}$. It then directly follows from Theorem~\ref{thm:embedding} that the origin of the nonlinear system \eqref{eqn:CL_NLsys_LPV} (which is equivalent to \eqref{eqn:CL_NLsys_data} or \eqref{eqn:CL_NLsys} on $\mathcal{X}$) is asymptotically stable. To this end, notice from \eqref{eqn:defGi} that
	\begin{equation}
		R_i = G_i \begin{bmatrix}
			M & 0\\ C_1 & C_2
		\end{bmatrix}.\label{eqn:Ri}
	\end{equation}
	Substituting this into the top left block of \eqref{eqn:P1a} gives 
	\[
	\sym\left( \begin{bmatrix}
		X_1\\ \mathcal{L}_i Z_0
	\end{bmatrix} G_i \begin{bmatrix}
		M & 0 \\ C_1 & C_2
	\end{bmatrix}\begin{bmatrix}
		\delta I_n & 0 \\ 0 & I_{S-n}
	\end{bmatrix} \right) + \begin{bmatrix}
		Q_j & 0\\ 0 & 0
	\end{bmatrix}.
	\]
	Letting $F\coloneqq\delta M,\,C_3 \coloneqq \delta C_1$ and recalling that $Z_0G_i=I_{S}$ (see \eqref{eqn:enforced_Z0Gi}), one can write
	\[
	\sym\left( \begin{bmatrix}
		X_1 G_i\\ \mathcal{L}_i
	\end{bmatrix} \begin{bmatrix}
		F & 0 \\ C_3 & C_2
	\end{bmatrix}\right) + \begin{bmatrix}
		Q_j & 0\\ 0 & 0
	\end{bmatrix},
	\]
	which is the top left block of \eqref{eqn:polyquadadmissibility_proof}. Similarly, substituting \eqref{eqn:Ri} in the off-diagonal blocks of \eqref{eqn:P1a} results in
	\[
	\begin{bmatrix}
		X_1 \\ \mathcal{L}_i Z_0
	\end{bmatrix}G_i\begin{bmatrix}
		M & 0 \\ C_1 & C_2
	\end{bmatrix}\begin{bmatrix}
		I_n \\ 0
	\end{bmatrix} + \begin{bmatrix}
		\delta M^\top\\ 0
	\end{bmatrix}.
	\]
	Again, noting that $Z_0G_i=I_{S}$ and recalling that $F=\delta M$ from before, we obtain the off-diagonal blocks of \eqref{eqn:polyquadadmissibility_proof}. As a result, \eqref{eqn:polyquadadmissibility_proof} holds and, hence, the system \eqref{eqn:embedded_CLsys} is poly-quadratically admissible (see Theorem~\ref{thm:LPV_polyquadadmissibility}). The claim now follows from Theorem~\ref{thm:embedding}.
\end{proof}

Theorem~\ref{thm:Kx} provides sufficient LMI conditions such that a polytopic stabilizing controller $u=K(x)Z(x)$ for the nonlinear system \eqref{eqn:NLsys} can be designed from data. Compared to existing works on data-driven stabilization of nonlinear systems, we make the following remarks. First, the majority of existing results yield a controller of the form $u=KZ(x)$, compare \cite{martin2023guarantees,depersis2023annrev,hu2024enforcing} and references therein, i.e., with a constant gain matrix $K$. In contrast, the proposed descriptor embedding allows us to design controllers of the form $u=K(x)Z(x)$, i.e., ones in which the gain matrix $K(x)$ is state-dependent. Although some works have featured such control gains, e.g., for stabilization of polynomial systems \cite{Guo2022polynomial}, obtaining the gain matrix there required solving a sum-of-squares program, which can become computationally challenging for higher dimensional systems. Theorem~\ref{thm:Kx} only requires solving a set of data-dependent LMIs that can be easily solved and returns (potentially globally) stabilizing controllers for systems that are \emph{not} restricted to polynomial ones (see Example~\ref{ex:global_inv_pend}).

Furthermore, many existing results on data-driven stabilization rely on minimizing and/or dominating the effect of the nonlinear terms. Despite being a simple and powerful technique, in many cases it might be beneficial to make use of the existing nonlinearities rather than minimizing or dominating their effect. This is especially true in the presence of uncertainties or when global nonlinearity cancellation is not possible (see Examples~\ref{ex:ex_nominal}-\ref{ex:ex_robust} below, where the proposed approach results in closed-loop systems with a larger region of attraction compared to existing works).

Before illustrating the results with examples, we comment on the complexity of the convex program \eqref{eqn:P1}. Notice that the number of LMIs that need to be solved scales quadratically with the number of vertices $N$ from (A2.3). Specifically, \eqref{eqn:P1} includes $N^2+N$ constraints where $N^2$ of them is due to \eqref{eqn:P1a} while the remaining $N$ are due to \eqref{eqn:P1b}. The reason for this is the use of a Lyapunov function whose parameter matrix $P(x)$ takes the form of a convex combination of positive definite matrices $Q^{-1}_i$, each of which corresponds to some vertex of the closed-loop polytope. To reduce the computational burden, one can use a uniform positive definite matrix $Q$ for all vertices and obtain a quadratic Lyapunov function of the form $V(x)=x^\top P x$, where $P=Q^{-1}$. This significantly reduces the number of LMIs required to solve (from $N^2+N$ to $2N$). To arrive at such a result, an extension of Theorem~\ref{thm:LPV_polyquadadmissibility} to the case of a uniform $Q$ matrix is used. Such a result was not reported in \cite{bara2011dilated}, but showing it is straightforward and follows similar arguments as in the proof of \cite[Th. 4.3(ii)]{bara2011dilated} (and is therefore omitted here). Based on this, we have the following corollary of Theorem~\ref{thm:Kx}.
\begin{corollary}\label{cor:Kx}
	Consider a nonlinear system as in \eqref{eqn:NLsys} and let Assumptions~\ref{asmp:f}-\ref{asmp:L} hold. Fix some $\delta\in\mathbb{R}$ and consider the following LMIs in the decision variables $R_i\in\mathbb{R}^{T\times S}$ for all $i\in\mathcal{I}=\{1,\ldots,N\}$, and $Q=Q^\top\succ0\in\mathbb{R}^{n\times n},\, M\in\mathbb{R}^{n\times n},\,C_1\in\mathbb{R}^{S-n \times n}$ and an invertible $C_2\in\mathbb{R}^{S-n \times S-n}$
	\begin{align}
		&\resizebox{\columnwidth}{!}{$
			\begin{bmatrix}
				\sym\left(\begin{bmatrix}
					X_1\\ \mathcal{L}_i Z_0
				\end{bmatrix}R_i\begin{bmatrix}
					\delta I_{n} & 0\\0 & I_{S-n}
				\end{bmatrix}\right) + \begin{bmatrix}
					Q & 0\\0 & 0
				\end{bmatrix} & \star \\ 
				\left(\begin{bmatrix}
					X_1\\ \mathcal{L}_i Z_0
				\end{bmatrix}R_i\begin{bmatrix}
					I_{n}\\ 0
				\end{bmatrix}+\begin{bmatrix}
					\delta M^\top\\ 0
				\end{bmatrix}\right)^\top & \sym(M)-Q
			\end{bmatrix}\succ0,
			$}\notag\\
		&Z_0 R_i = \begin{bmatrix}
			M & 0 \\ C_1 & C_2
		\end{bmatrix}. \label{eqn:corP1}
	\end{align}
	If \eqref{eqn:corP1} is feasible, then the origin of \eqref{eqn:CL_NLsys}, with $K(x)=\sum_{i=1}^{N}\alpha_i(x)K_i$ and $K_i$ as in \eqref{eqn:P1_gain}, is locally asymptotically stable. A Lyapunov function for the closed-loop system is given by $V(x_t)=x_t^\top Q^{-1} x_t$. Furthermore, if $\mathcal{X}=\mathbb{R}^n$, then the origin is globally asymptotically stable. \hfill$\square$
\end{corollary}

Both Theorem~\ref{thm:Kx} and Corollary~\ref{cor:Kx} provide sufficient conditions for obtaining stabilizing controllers for the nonlinear system that generated the data. The difference is that the latter certifies stability of the closed-loop system with a uniform quadratic Lyapunov function. This may result in different ROA estimates when the resulting controller is only locally stabilizing. Notice that if a uniform $Q$ exists satisfying \eqref{eqn:corP1}, then \eqref{eqn:P1} is satisfied with $Q_i=Q$ for all $i$. Conversely, one can easily show that if \eqref{eqn:P1} is feasible, then $Q=\frac{1}{N}\sum_{i=1}^N Q_i$ satisfies~\eqref{eqn:corP1}.

\subsection{Global vs local stabilization}
\label{sec:global_vs_local}
The controller obtained from Theorem~\ref{thm:Kx} (or Corollary~\ref{cor:Kx}) is globally stabilizing if Assumption~\ref{asmp:L} holds globally, i.e., with $\mathcal{X}=\mathbb{R}^n$. The following example illustrates the results of Theorem~\ref{thm:Kx} in such a setting. Unlike the approach from \cite{DePersis2023}, the resulting controller does \textit{not} result in nonlinearity cancellation but still globally asymptotically stabilizes the~origin.

\begin{example}\label{ex:global_inv_pend}
	Consider a discretized model of an inverted pendulum
	\[
	\begin{aligned}
		x_{1,t+1} &= x_{1,t} + T_s x_{2,t}\\
		x_{2,t+1} &= \left(1-\frac{T_s\mu}{\bar{m}\ell^2}\right)x_{2,t} + \frac{T_s g}{\ell}\sin(x_{1,t}) + \frac{T_s}{\bar{m}\ell^2}u_t,
	\end{aligned}
	\]
	where $x_{1,t},\,x_{2,t}$ are the angular position and velocity of the pendulum, respectively, $T_s=0.1$ is the sampling time, $\bar{m}=1$ is the mass of the pendulum, $\ell=1$ is the length of the pendulum, $\mu=0.01$ is the coefficient of static friction, and $g=9.81$ is the gravitational acceleration. We choose as basis functions $Z(x)=\begin{bmatrix}
		x_{1} & x_{2} & \sin(x_1)
	\end{bmatrix}^\top$ which satisfy Assumption~\ref{asmp:f}. For such a choice of $Z$, the following is a function $L$ satisfying Assumption~\ref{asmp:L}
	\[
	L(x_1) = \begin{bmatrix}
		\mathrm{sinc}(x_1) & 0 & -1
	\end{bmatrix},
	\]
	where
	\[
	\mathrm{sinc}(x_1) = \begin{cases}
		\sin(x_1)/x_1, \quad & x_1\neq0,\\
		1, \quad & x_1 = 0.
	\end{cases}
	\]
	It is easy to see that $L$ satisfies (A2.1) and (A2.2). To see that (A2.3) holds, notice that we can express $L$ as $L(x) = \alpha_1(x) \mathcal{L}_1 + \alpha_2(x) \mathcal{L}_2$, for all $x\in\mathbb{R}^n$, where
	\[
	\begin{aligned}
		\mathcal{L}_1 &= L(v_1),\\
		\alpha_1(x) &= \frac{\mathrm{sinc}(x_1) - \mathrm{sinc}(v_2)}{\mathrm{sinc}(v_1) - \mathrm{sinc}(v_2)},
	\end{aligned}\,\,\,
	\begin{aligned}
		\mathcal{L}_2 &= L(v_2),\\
		\alpha_2(x) &= \frac{\mathrm{sinc}(v_1) - \mathrm{sinc}(x_1)}{\mathrm{sinc}(v_1) - \mathrm{sinc}(v_2)},
	\end{aligned}
	\]
	with $v_1=0$ and $v_2 = 4.4934$ being the points at which $\mathrm{sinc}$ attains its maximum ($\mathrm{sinc}(v_1)=1$) and minimum ($\mathrm{sinc}(v_2)=-0.2172$) values, respectively. Notice that here, Assumption~\ref{asmp:L} holds globally with $N=2$ vertices, with the functions $\alpha_i$ obtained by noting that $\mathrm{sinc}$ is a bounded function that can be written in terms of its maximum and minimum values.
	
	To apply the results of Theorem~\ref{thm:Kx}, we run an offline experiment and collect $\{u_k^d,x_k^d\}_{k=0}^{T}$ of length $T=10$, where the input is sampled uniformly from $[-0.5,0.5]$ and initial conditions for each state sampled from the same interval. We set $\delta = 0.5$ and find that \eqref{eqn:P1} is feasible and returns the following controller
	\begin{equation}
		\begin{aligned}
			K(x) &= \alpha_1(x)K_1 + \alpha_2(x)K_2,\\
			K_1 &= \begin{bmatrix}
				-8.4519 & -11.5047 & -10.0532
			\end{bmatrix},\\
			K_2 &= \begin{bmatrix}
				-8.6412 & -11.4817 &  -9.5221
			\end{bmatrix}.
		\end{aligned}
		\label{eqn:ex1_global_ourK}
	\end{equation}
	Comparing this with the controller obtained from the nonlinearity minimization approach \cite{DePersis2023}
	\begin{equation}
		K_{\mathrm{NLmin}} = \begin{bmatrix}
			-5.6141 & -1.5959 & -9.8100
		\end{bmatrix}
		\label{eqn:ex1_global_NLminK}
	\end{equation}
	we clearly see that the latter performs cancellation (as evident by the last entry which cancels the effect of $\frac{T_s g}{\ell}\sin(x_{1,t})$ in the dynamics of $x_2$). In contrast, the controller \eqref{eqn:ex1_global_ourK} obtained from Theorem~\ref{thm:Kx} does not result in nonlinearity cancellation, but still globally asymptotically stabilizes the origin of the nonlinear system. This can be certified by the following Lyapunov function, which is also obtained from \eqref{eqn:P1}
	\begin{align}
		V(x) &= x^\top P(x) x, \quad P(x)=\sum\nolimits_{i=1}^{2}\alpha_i(x)P_i,\notag\\
		P_1 &= Q_1^{-1} = \begin{bmatrix}
			0.1793  &  0.0536\\
			0.0536  &  0.0611
		\end{bmatrix},\notag\\
		P_2 &= Q_2^{-1} = \begin{bmatrix}
			0.1793  &  0.0537\\
			0.0537  &  0.0611
		\end{bmatrix}.\tag*{$\square$}
	\end{align}
\end{example}

We highlight the fact that the resulting controller is \textit{not} a canceling controller. Instead, it exploits the inherent nonlinearities of the system to later successfully globally stabilize the origin. Global stabilization without cancellation can be achieved for certain classes of nonlinear systems, e.g., polynomial systems. However, to the best of our knowledge, no other technique exists for global stabilization of systems satisfying Assumptions~\ref{asmp:f} and \ref{asmp:L} which does not resort to cancellation of the nonlinearities (as the one considered in Example~\ref{ex:global_inv_pend}).

For the previous example, global stabilization was possible since Assumption~\ref{asmp:L} was satisfied with $\mathcal{X}=\mathbb{R}^n$. In Section~\ref{sec:embedding}, we discussed Assumption~\ref{asmp:L} in detail and showed in Theorem~\ref{thm:asmpLmonomoialZ} that Assumption~\ref{asmp:L} is satisfied locally when choosing $Z$ as all monomials of the state up to some finite degree. In that case, the set $\mathcal{X}\subset\mathbb{R}^n$ is a polytope and the resulting controller obtained using Theorem~\ref{thm:Kx} is no longer guaranteed to be globally stabilizing. Nevertheless, one can potentially choose $v_i\in\mathcal{X}$ in Assumption~\ref{asmp:L} such that it covers a practical operating range of the system.

In what follows, we investigate the properties of the obtained local controllers. In particular, we illustrate how an estimate of the region of attraction of the closed-loop system can be obtained using only available data and the outcome of Theorem~\ref{thm:Kx} (or Corollary~\ref{cor:Kx}). To this end, suppose \eqref{eqn:P1} is feasible and consider the Lyapunov function $V(x_t)=x_t^\top P(x_t) x_t$ along with its difference
\begin{align}
    &V(x_{t+1}) - V(x_t)\notag\\
    &= x_{t+1}^\top P(x_{t+1}) x_{t+1} - x_t^\top P(x_t) x_t\notag\\
    &\stackrel{\eqref{eqn:CL_NLsys_LPV}}{=} \left(X_1G(x_t)Z(x_t)\right)^\top P(x_{t+1}) X_1G(x_t)Z(x_t) - x_t^\top P(x_t) x_t\notag\\
    &\eqqcolon h(x_t),\label{eqn:Vdiff_data}
\end{align}
where 
\[
\begin{aligned}
    P(x_{t+1}) &= \sum\nolimits_{i=1}^{N}\alpha_i(x_{t+1}) Q_i^{-1}\\
    &= \sum\nolimits_{i=1}^{N}\alpha_i\left(X_1G(x_t)Z(x_t)\right) Q_i^{-1}.
\end{aligned}
\]
Notice that all these quantities are known and, hence, the function $h(x)$ is also known. Let $\mathcal{V}\coloneqq\left\lbrace x ~|~ h(x) < 0 \right\rbrace\cup\{0\}$ denote the set on which the Lyapunov function difference is negative definite and let $\mathcal{R}_{c}=\left\lbrace x ~|~ V(x)\leq c\right\rbrace$ denote a sublevel set of the Lyapunov function. It is well known that any sublevel set contained in $\mathcal{V}\cap\mathcal{X}$ is a positively invariant set and represents an estimate of the ROA of the closed-loop system, compare~\cite{DePersis2023,haddad2008nonlinear}. Depending on the size of the set $\mathcal{X}$, the above procedure may yield small ROA estimates. This, however, does not necessarily mean that the true ROA is small. This is because the intersection with $\mathcal{X}$ is only needed for analysis purposes and for obtaining a guaranteed ROA estimate, whereas the closed-loop system (and in particular the gain $K(x)$) is still well-defined outside~$\mathcal{X}$.

In the following, we illustrate the proposed methods from Theorem~\ref{thm:Kx} and Corollary~\ref{cor:Kx} when only local controllers are obtained, and compare them to the nonlinearity minimization approach from \cite{DePersis2023}. 

\begin{example}
\label{ex:ex_nominal}
Consider the following polynomial system
\[
\begin{aligned}
    x_{1,t+1} &= 0.8 x_{2,t} + 0.2 x_{1,t}^3,\\
    x_{2,t+1} &= -0.6 x_{1,t} + x_{2,t}^2 - u_t.
\end{aligned}
\]
For such a system, exact nonlinearity cancellation is not possible, due to the $0.2x_1^3$ term in the dynamics of the first state. In this example, we will compare three different controllers: two polytopic controllers of the form $u=K(x)Z(x)$ using (i) Theorem~\ref{thm:Kx}, (ii) Corollary~\ref{cor:Kx}, and (iii) a controller of the form $u=K_{\mathrm{NLmin}}Z(x)$ using the approach from \cite{DePersis2023}.

To begin, consider the following choice of basis functions $Z$ which satisfy Assumption~\ref{asmp:f}
\[
Z(x) = \begin{bmatrix}
    x_1 & x_2 & x_1^2 & x_2^2 & x_1^3 & x_2^3
\end{bmatrix}^\top.
\]
For such a choice of $Z$, we can choose $L$ as
\[
L(x) = \begin{bmatrix}
    x_1 & 0 & -1 & 0 & 0 & 0\\
    0 & x_2 & 0 & -1 & 0 & 0\\
    0 & 0 & x_1 & 0 & -1 & 0\\
    0 & 0 & 0 & x_2 & 0 & -1\\
\end{bmatrix}
\]
which satisfies Assumption~\ref{asmp:L} only locally (compare\footnote{The function $Z$ chosen in this example does not contain all momonials up to degree 3; nevertheless, the function $L$ can still be determined using the same procedure as in the proof of Theorem~\ref{thm:asmpLmonomoialZ}.} Theorem~\ref{thm:asmpLmonomoialZ}). In particular, for $\mathcal{X} = \lbrace x\in\mathbb{R}^2 ~|~ |x_i| \leq \gamma_i,\, \gamma_i\geq0 \rbrace\subset\mathbb{R}^2$, Assumption (A2.3) holds with $N=4$ vertices 
\[
\begin{aligned}
    \mathcal{L}_1 &= L([\gamma_1 \,\, \gamma_2]^\top),\\
    \mathcal{L}_3 &= L([\gamma_1 \,\, -\gamma_2]^\top),
\end{aligned}\quad
\begin{aligned}
    \mathcal{L}_2 &= L([-\gamma_1 \,\, \gamma_2]^\top),\\
    \mathcal{L}_4 &= L([-\gamma_1 \,\, -\gamma_2]^\top),
\end{aligned}
\]
and weights $\alpha_i(x)$ obtained by the generalized barycentric coordinates of the polytope $\mathcal{X}$ (see \cite{floater2015generalized} and Appendix~\ref{app:exmpL})
\[
\begin{aligned}
    \alpha_1(x) &= \frac{(\gamma_1+x_1)(\gamma_2+x_2)}{4\gamma_1\gamma_2},\\
    \alpha_3(x) &= \frac{(\gamma_1+x_1)(\gamma_2-x_2)}{4\gamma_1\gamma_2},
\end{aligned}\,\,
\begin{aligned}
    \alpha_2(x) &= \frac{(\gamma_1-x_1)(\gamma_2+x_2)}{4\gamma_1\gamma_2},\\
    \alpha_4(x) &= \frac{(\gamma_1-x_1)(\gamma_2-x_2)}{4\gamma_1\gamma_2}.
\end{aligned}
\]

All three controllers are designed using the same offline collected data set $\{u_k^d,x_k^d\}_{k=0}^{T}$ of length $T=10$, where the input is sampled uniformly from $[-0.5,0.5]$ and initial conditions for each state sampled from the same interval. For the first two controllers, we set $\gamma_1=3,\gamma_2=4$ and set\footnote{The LMIs \eqref{eqn:P1} (or \eqref{eqn:corP1}) were also feasible for other values of $\delta$. Here, as well as in the following examples, we report the values of $\delta$ that resulted in the largest ROA estimates. In particular, for each fixed value of $\delta\in\{10^k ~|~ k \text{ is an integer in } [-6,4]\}\cup\{0\}$, we solve a semi-definite program whose objective function incentivizes larger ROA estimates, subject to \eqref{eqn:P1} (or \eqref{eqn:corP1}). Specifically, we aimed to maximize the volume of the ellipsoid $x^\top (\sum_{i=1}^{N}\alpha_i(x)Q_i^{-1}) x$ (or $x^\top Q^{-1} x$) by minimizing $-\sum_{i=1}^{N}\mathrm{log}\,\mathrm{det}(Q_i)$ (or $-\mathrm{log}\,\mathrm{det}(Q)$), compare~\cite{Boyd2004}. For the nonlinearity minimization approach, we solve \cite[Eq. (30)]{DePersis2023} which also incentivizes larger ROA estimates.} $\delta=0$ in both \eqref{eqn:P1} and \eqref{eqn:corP1}. All programs were found feasible and return stabilizing controllers. 
\begin{figure}[t!]
    \centering
    \includegraphics[width=0.9\linewidth]{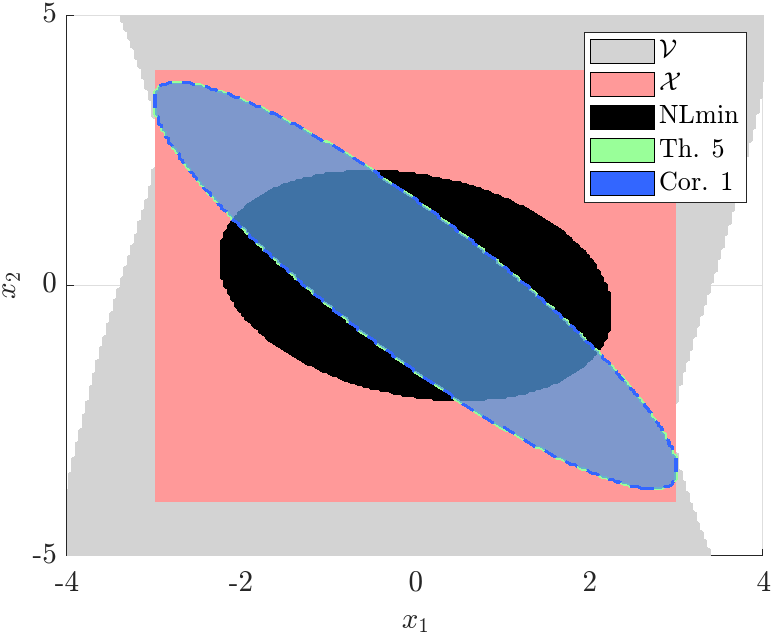}
    \caption{ROA estimates using all three methods (NLmin refers to the method from \cite{DePersis2023}). For each method, we show the largest sub level set $\mathcal{R}_{c_{\max}}$ contained in its corresponding set $\mathcal{V}$ (the latter shown only for the methods proposed in the paper and not for the NLmin approach). The areas of each region are as follows: NLmin: 14.6960 (in black), Theorem~\ref{thm:Kx}: 15.4224 (in green) representing approximately 5\% increase over NLmin, and Corollary~\ref{cor:Kx}: 15.4224 (in blue) also representing approximately 5\% increase over NLmin.}
    \label{fig:ROAoverlayed}
\end{figure}
\begin{figure}[t!]
	\centering
	\includegraphics[width=0.9\linewidth]{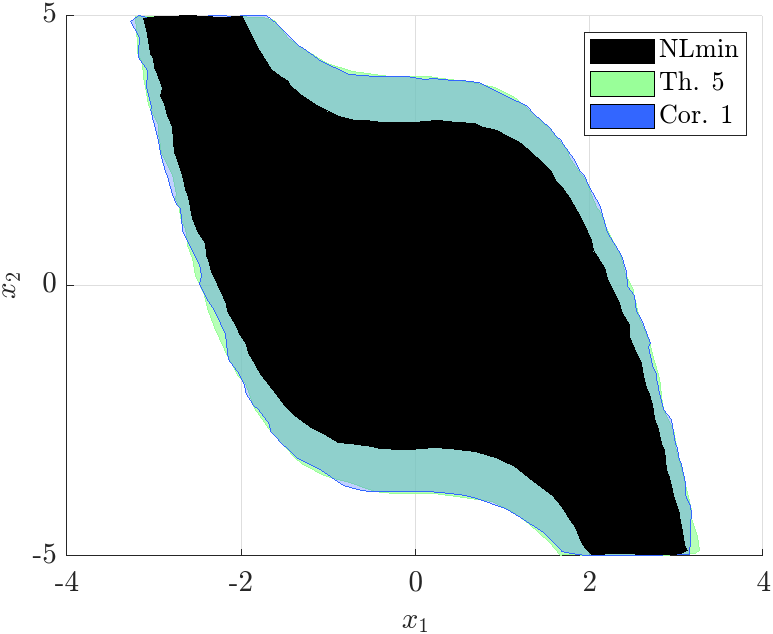}
	\caption{Numerical estimates of the true ROA for all three methods. The areas of each region are as follows: NLmin: 31.8450, Theorem~\ref{thm:Kx}: 40.2670 representing approximately 26.4\% increase over NLmin, and Corollary~\ref{cor:Kx}: 40.1383, approximately a 26\% increase over NLmin.}
	\label{fig:numericalROA}
\end{figure}

Figure~\ref{fig:ROAoverlayed} shows the data-based ROA estimates for the closed-loop system under each controller. It can be seen that the controllers designed using the proposed framework (i.e., those from Theorem~\ref{thm:Kx} and Corollary~\ref{cor:Kx}) both result in a larger ROA estimate than the method from \cite{DePersis2023}. Next, we numerically estimate the true ROA resulting from application of each controller. To do this, we run 5000 closed-loop experiments each of which is 200 steps long and record the initial conditions for which the experiment has converged. Convergence here is determined by testing whether the values of the state reached machine precision, i.e., $\leq 10^{-6}$ in 200 steps. It can be seen in Figure~\ref{fig:numericalROA} that the proposed controllers enjoy a larger ROA compared to the approach from~\cite{DePersis2023}.\hfill$\square$
\end{example}

In the following section, we consider various extensions to the proposed approach. In particular, we consider the case of noisy data and inexact basis function decomposition, then briefly discuss a more general class of nonlinear systems.
\vfill
\section{Extensions}\label{sec:extensions}
\subsection{Inexact basis function expansion and robust stabilization}\label{sec:inexact}
Inevitably, measured data is affected by noise and, in practice, one may not necessarily have sufficient knowledge about the physics of the system to properly choose basis functions $Z$ that satisfy Assumption~\ref{asmp:f}. For the latter, one can still write
\[
f(x) = \widehat{A}Z(x) + \varepsilon^f(x),
\]
where $\varepsilon^{f}$ is the approximation error and $\widehat{A}$ is an \textit{unknown} matrix that minimizes the average approximation error on a (in general compact) subset $\mathcal{S}\subset\mathbb{R}^n$ (with $0\in\mathrm{int}(\mathcal{S})$)
\[
\widehat{A} \coloneqq \argmin_A \int_{\mathcal{S}} \norm{f(x)-AZ(x)}^2 dx.
\]
We assume that $\varepsilon^f$ is uniformly bounded on $\mathcal{S}$, i.e.,
\begin{equation}
    \norm{\varepsilon^f(x)} \leq \bar\varepsilon,\quad\forall x\in\mathcal{S},\label{eqn:pointwisebnd_epsilon}
\end{equation}
for some $\bar\varepsilon>0$. This is the case, e.g., if $\mathcal{S}$ is compact (since $f$ and $Z$ are assumed to be continuous). We formalize this in the following assumption (which relaxes Assumption~\ref{asmp:f}).
~
\begin{assumption}\label{asmp:f2}
    The function $f$ can be written as $f(x)=\widehat{A}Z(x)+\varepsilon^f(x)$, where $\widehat{A}\in\mathbb{R}^{n\times S}$ is an unknown matrix, $Z:\mathbb{R}^n\to\mathbb{R}^{S}$ is a known, continuous, vector-valued function that takes the form $Z(x) = \begin{bmatrix}
        x^\top & \underline{Z}(x)^\top
    \end{bmatrix}^\top$, with $Z(0)=0$, and
    the approximation error $\varepsilon^f$ is uniformly bounded, i.e., there exists a known $\bar\varepsilon>0$ such that \eqref{eqn:pointwisebnd_epsilon} holds.\hfill$\square$
\end{assumption}
~
A poor choice of basis functions can obviously result in a large $\varepsilon^f$. Alternatively, choosing $Z$ as all monomials up to some finite degree results in a tight approximation of arbitrary smooth functions $f$ on compact sets \cite{Stone1948} and, hence, yields a sufficiently small $\varepsilon^f$ that can be properly accounted for. Moreover, such a choice of basis functions is advantageous since, for such $Z$ one can always find a function $L$ that satisfies Assumption~\ref{asmp:L} as was shown in Theorem~\ref{thm:asmpLmonomoialZ}.

The goal of this section is to design a polytopic state-feedback controller $u=K(x)Z(x)$, with $K(x)$ as in \eqref{eqn:CLgain}, such that the origin of the closed-loop system 
\begin{equation}
	x_{t+1} = (\widehat{A}+BK(x_t))Z(x_t) + \varepsilon^f(x_t) \label{eqn:robust_CL_NLsys}
\end{equation}
is locally asymptotically stable. Since $\widehat{A},B$ are unknown, we assume availability of \textit{noisy} data collected from an open-loop experiment $\{u_k^d,\,y_k^d\}_{k=0}^{T}$ where $y_k^d = x_k^d + w_k^d$. We make the following assumption on the noise sequence.
\begin{assumption}\label{asmp:noise}
	There exists a known $\bar w$ such that  $\norm{w_k}\leq\bar w$ for all $k\in\mathbb{N}$.\hfill$\square$
\end{assumption}

Let the collected data be arranged as follows
\begin{subequations}
  \begin{align}
   &\begin{aligned}
       U_0 &= \begin{bmatrix}
            u_0^d & u_1^d & \cdots & u_{T-1}^d
        \end{bmatrix}\in\mathbb{R}^{m\times T},\\
        Z_0 &= \begin{bmatrix}
            Z(y_0^d) & Z(y_1^d) & \cdots & Z(y_{T-1}^d)
        \end{bmatrix}\in\mathbb{R}^{S\times T},\\
        Y_1 &= \begin{bmatrix}
            y_1^d & y_2^d & \cdots & y_{T}^d
        \end{bmatrix}\in\mathbb{R}^{n\times T}, \\
   \end{aligned}
 \intertext{and consider the corresponding \textit{unknown} matrices}
        &\begin{aligned}
            \widetilde{Z}_0 &= \begin{bmatrix}
            Z(x_0^d) & Z(x_1^d) & \cdots & Z(x_{T-1}^d)
        \end{bmatrix}\in\mathbb{R}^{S\times T},\\
        \mathcal{E}_0 &= \begin{bmatrix}
    \varepsilon^f(x_0^d) & \varepsilon^f(x_1^d) & \cdots & \varepsilon^f(x_{T-1}^d)
\end{bmatrix}\in\mathbb{R}^{n\times T},\\
        W_1 &= \begin{bmatrix}
            w_1^d & w_2^d & \cdots & w_{T}^d
        \end{bmatrix}\in\mathbb{R}^{n\times T}.
        \end{aligned}
    \end{align}%
    \label{eqn:noisy_data_mats}%
\end{subequations}%
By straightforward manipulations, one can see that the data matrices satisfy
\begin{equation}
    Y_1 = \widehat{A}Z_0 + BU_0 + D_0,
    \label{eqn:noisy_state_identity}
\end{equation}
where $D_0\coloneqq W_1 + \mathcal{E}_0 + \widehat{A}(\widetilde{Z}_0-Z_0)$. Due to Assumptions~\ref{asmp:f2} and \ref{asmp:noise}, the matrix $D_0$ belongs to a set $\mathcal{D}\coloneqq\{D\in\mathbb{R}^{n\times T}~|~DD^\top\preceq\Delta\}$ for some $\Delta\succeq0\in\mathbb{R}^{n\times n}$ which is assumed to be known. The following is an analogous result to Lemma~\ref{lemma:DBrep} and, therefore, its proof is omitted for brevity. 
\begin{lemma}\label{lemma:robust_DBrep}
	Consider matrices $K_i\in\mathbb{R}^{m\times S}, G_i\in\mathbb{R}^{T\times S}$, for $i\in\mathcal{I}\coloneqq\{1,\ldots,N\}$, satisfying \eqref{eqn:KGrelation}. Then, \eqref{eqn:robust_CL_NLsys} can be equivalently written as
	\begin{equation}
		x_{t+1} = (Y_1-D_0)G(x_t)Z(x_t) + \varepsilon^f(x_t) \label{eqn:robust_CL_NLsys_data}
	\end{equation}
	where $G(x) = \sum\nolimits_{i=1}^{N}\alpha_i(x)G_i$ and $\alpha_i$ defined below \eqref{eqn:CLgain}.\hfill$\square$
\end{lemma}

Notice that the data-based representation \eqref{eqn:robust_CL_NLsys_data} accounts for two sources of uncertainty. Namely, (i) neglected nonlinearities and noisy offline data captured by $D_0$, and (ii) neglected nonlinearities in closed-loop operation captured by $\varepsilon^f(x_t)$. Analogously to Theorem~\ref{thm:embedding}, it can be shown that poly-quadratic admissibility of a polytopic LPVd system of the form
\begin{equation}
	\begin{bmatrix}
		I_n & 0\\ 0 & 0
	\end{bmatrix}\eta_{t+1}= \begin{bmatrix}
		(Y_1-D_0)G(\rho_t) \\ L(\rho_t)
	\end{bmatrix}\eta_t,
	\label{eqn:robust_embedded_CLsys}
\end{equation}
for which $\rho_t\in\mathcal{X}\subseteq\mathbb{R}^n$, implies asymptotic stability of the origin of
\begin{equation}
	\begin{aligned}
		x_{t+1} &= (Y_1-D_0)G(x_t)Z(x_t)\\
		&= (Y_1-D_0)\left(\overline{G}(x_t) - \underline{G}(x_t)L_2^{-1}(x_t)L_1(x_t)\right)x_t
	\end{aligned}\label{eqn:robust_linearsubsystem}
\end{equation}
where the second equality holds due to (A2.1) and (A2.2), i.e., that $\underline{Z}(x_t) = -L_2^{-1}(x)L_1(x)$, with $G(x) = \begin{bmatrix} \overline{G}(x) & \underline{G}(x)\end{bmatrix}$ partitioned appropriately (compare \eqref{eqn:CL_NLsys_LPV}). Unlike Theorem~\ref{thm:embedding}, however, this alone is not sufficient to conclude asymptotic stability of the origin of \eqref{eqn:robust_CL_NLsys_data} due to the presence of the neglected nonlinearities in closed-loop operation. However, under the additional assumption that $\varepsilon^f(x_t)$ goes to zero faster than linearly, i.e.,
\begin{equation}
	\lim\limits_{\norm{x}\to0}\frac{\norm*{\varepsilon^f(x)}}{\norm{x}} = 0,\label{eqn:asmp_epsilon}
\end{equation}
local asymptotic stability of the origin of \eqref{eqn:robust_CL_NLsys_data} (which is equivalent to \eqref{eqn:robust_CL_NLsys} on $\mathcal{X}$) directly follows. This is summarized in the following theorem, which represents a robust extension of the results of Theorem~\ref{thm:Kx} to account for measurement noise and neglected nonlinearities.
~
\begin{theorem}\label{thm:Kx_robust}
    Consider a nonlinear system as in \eqref{eqn:NLsys} and let Assumptions~\ref{asmp:L}-\ref{asmp:noise} and \eqref{eqn:asmp_epsilon} hold. Fix some scalar $\delta\in\mathbb{R}$ and consider \eqref{eqn:P2} (shown at the top of the next page) in the decision variables $\lambda_i>0,\,R_i\in\mathbb{R}^{T\times S},\,Q_i=Q_i^\top\succ0\in\mathbb{R}^{n\times n}$, for all $i\in\mathcal{I}=\{1,\ldots,N\}$, $M\in\mathbb{R}^{n\times n},\,C_1\in\mathbb{R}^{S-n \times n}$ and an invertible $C_2\in\mathbb{R}^{S-n \times S-n}$.
    \begin{figure*}[!t]
        \centering
        \begin{subequations}
            \begin{align}
            &\begin{bmatrix}
                \sym\left(\begin{bmatrix}
                    Y_1\\ \mathcal{L}_i Z_0
                \end{bmatrix}R_i\begin{bmatrix}
                    \delta I_{n} & 0\\0 & I_{S-n}
                \end{bmatrix}\right) + \begin{bmatrix}
                    Q_j - \lambda_i \Delta & 0\\0 & 0
                \end{bmatrix} & \star & \star\\ 
                \left(\begin{bmatrix}
                    Y_1\\ \mathcal{L}_i Z_0
                \end{bmatrix}R_i\begin{bmatrix}
                    I_{n}\\ 0
                \end{bmatrix}+\begin{bmatrix}
                    \delta M^\top\\ 0
                \end{bmatrix}\right)^\top & \sym(M)-Q_i & \star\\
                R_i\begin{bmatrix}
                    \delta I_{n} & 0\\ 0 & I_{S-n}
                \end{bmatrix} & R_i\begin{bmatrix}
                    I_n\\0
                \end{bmatrix} & \lambda_i I_{T}
            \end{bmatrix}\succ0,\quad\forall (i,j)\in\mathcal{I}\times\mathcal{I},\label{eqn:P2a}\\
            &Z_0 R_i = \begin{bmatrix}
                M & 0 \\ C_1 & C_2
            \end{bmatrix}.\label{eqn:P2b}
        \end{align}\label{eqn:P2}
        \end{subequations}
        \hrule
    \end{figure*}
    If \eqref{eqn:P2} is feasible, then the origin of \eqref{eqn:robust_CL_NLsys}, with $K(x)=\sum_{i=1}^{N}\alpha_i(x)K_i$ and
    \begin{equation}
        K_i = U_0 R_i \begin{bmatrix}
            M & 0 \\ C_1 & C_2
        \end{bmatrix}^{-1},\label{eqn:P2_gain}
    \end{equation}
    is locally asymptotically stable. A Lyapunov function for the closed-loop system is given by $V(x_t)=x_t^\top P(x_t)x_t$ with $P(x_t) = \sum_{i=1}^{N}\alpha_i(x_t)Q_i^{-1}$.\hfill$\square$
\end{theorem}
\begin{proof}
	Let \eqref{eqn:P2} hold. The central block of \eqref{eqn:P2a} implies $\sym(M)-Q_i\succ0$ for all $i$, and since $Q_i\succ0$ by assumption, it follows that $M$ is invertible. Furthermore, since $C_2$ is invertible by assumption, then one can define
	\begin{equation}
		G_i \coloneqq R_i\begin{bmatrix}
			M & 0\\ C_1 & C_2
		\end{bmatrix}^{-1}.\label{eqn:robust_defGi}
	\end{equation}
	Pre-multiplying both sides from the left by $Z_0$ results in 
	\[
	Z_0G_i = Z_0R_i\begin{bmatrix}
		M & 0\\ C_1 & C_2
	\end{bmatrix}^{-1}\stackrel{\eqref{eqn:P2b}}{=} I_{S},
	\]
	which, together with \eqref{eqn:P2_gain}, enforces \eqref{eqn:KGrelation} and, hence, the data-based representation of the closed-loop system \eqref{eqn:robust_CL_NLsys_data} holds (see Lemma~\ref{lemma:robust_DBrep}).
	
    Next, we will show that \eqref{eqn:P2a} is equivalent to the following LMIs being satisfied
    \begin{equation}
        \resizebox{\columnwidth}{!}{$\begin{bmatrix}
            \sym\left(\begin{bmatrix}
                Y_1-D_0\\ \mathcal{L}_i Z_0
            \end{bmatrix}R_i\begin{bmatrix}
                \delta I_n & 0 \\ 0 & I_{S-n}
            \end{bmatrix}
            \right) + \begin{bmatrix}
                Q_j & 0 \\ 0 & 0
            \end{bmatrix} & \star\\
            \left(\begin{bmatrix}
                Y_1-D_0\\ \mathcal{L}_i Z_0
            \end{bmatrix}R_i\begin{bmatrix}
                I_n\\0
            \end{bmatrix}+\begin{bmatrix}
                \delta M^\top \\ 0
            \end{bmatrix}
            \right)^\top & \sym(M)-Q_i
        \end{bmatrix}\succ0,$}\label{eqn:LMIrobust_proof}
    \end{equation}
    for all $(i,j)\in\mathcal{I}\times\mathcal{I}$. This is because, by Theorem~\ref{thm:LPV_polyquadadmissibility}, this in turn implies poly-quadratic admissibility of \eqref{eqn:robust_embedded_CLsys} (compare also the proof of Theorem~\ref{thm:Kx}). Via a Schur complement and straightforward manipulations of the resulting matrices, \eqref{eqn:P2a} can be equivalently written as
    \[
    \begin{aligned}
        &\resizebox{\columnwidth}{!}{$\begin{bmatrix}
        \sym\left(\begin{bmatrix}
                    Y_1\\ \mathcal{L}_i Z_0
                \end{bmatrix}R_i\begin{bmatrix}
                    \delta I_{n} & 0\\0 & I_{S-n}
                \end{bmatrix}\right) + \begin{bmatrix}
                    Q_j & 0\\0 & 0
                \end{bmatrix} & \star\\
                \left(\begin{bmatrix}
                    Y_1\\ \mathcal{L}_i Z_0
                \end{bmatrix}R_i\begin{bmatrix}
                    I_{n}\\ 0
                \end{bmatrix}+\begin{bmatrix}
                    \delta M^\top\\ 0
                \end{bmatrix}\right)^\top & \sym(M)-Q_i
    \end{bmatrix}$}\\
    &-\lambda_i^{-1} \begin{bmatrix}
            \delta I_n & 0\\ 0 & I_{S-n} \\ I_n & 0
        \end{bmatrix} R_i^\top R_i\begin{bmatrix}
            \delta I_n & 0 & I_n\\0 & I_{S-n}&0
        \end{bmatrix}\\
        &-\lambda_i\begin{bmatrix}
        I_n \\ 0 \\0 
    \end{bmatrix}\Delta\begin{bmatrix}
        I_n & 0 & 0
    \end{bmatrix}\succ0.
    \end{aligned} 
    \]
    Since $\lambda_i>0$ and $\Delta\succeq0$, we can use the strict Petersen's lemma (see \cite{petersen1987,bisoffi2022data}) and write the previous inequality equivalently~as
    \[
    \begin{aligned}
        &\resizebox{\columnwidth}{!}{$\begin{bmatrix}
        \sym\left(\begin{bmatrix}
                    Y_1\\ \mathcal{L}_i Z_0
                \end{bmatrix}R_i\begin{bmatrix}
                    \delta I_{n} & 0\\0 & I_{S-n}
                \end{bmatrix}\right) + \begin{bmatrix}
                    Q_j & 0\\0 & 0
                \end{bmatrix} & \star\\
                \left(\begin{bmatrix}
                    Y_1\\ \mathcal{L}_i Z_0
                \end{bmatrix}R_i\begin{bmatrix}
                    I_{n}\\ 0
                \end{bmatrix}+\begin{bmatrix}
                    \delta M^\top\\ 0
                \end{bmatrix}\right)^\top & \sym(M)-Q_i
    \end{bmatrix}$}\\
    &-\begin{bmatrix}
        I_n \\ 0 \\0 
    \end{bmatrix}D R_i \begin{bmatrix}
            \delta I_n & 0 & I_n\\0 & I_{S-n}&0
        \end{bmatrix}\\
        &-\begin{bmatrix}
            \delta I_n & 0\\ 0 & I_{S-n} \\ I_n & 0
        \end{bmatrix} R_i^\top D^\top\begin{bmatrix}
        I_n & 0 & 0
    \end{bmatrix}\succ0,
    \end{aligned}
    \]
    for all $(i,j)\in\mathcal{I}\times\mathcal{I}$ and all $D$ satisfying $DD^\top\preceq\Delta$ (i.e., for all $D\in\mathcal{D}$). Since $D_0\in\mathcal{D}$, then the above also holds for $D_0$. Collecting the terms results in \eqref{eqn:LMIrobust_proof}, i.e., that \eqref{eqn:robust_embedded_CLsys} is poly-quadratically admissible, as desired.
    
    Following similar arguments as in the proof of Theorem~\ref{thm:embedding}, it is straightforward to show that poly-quadratic admissibility of \eqref{eqn:robust_embedded_CLsys} implies local asymptotic stability of \eqref{eqn:robust_linearsubsystem}. This, together with \eqref{eqn:asmp_epsilon}, shows local asymptotic stability of the origin of \eqref{eqn:robust_CL_NLsys_data} (which is equivalent to \eqref{eqn:robust_CL_NLsys} on $\mathcal{X}$).
\end{proof}

Theorem~\ref{thm:Kx_robust} provides sufficient conditions for obtaining locally stabilizing controllers for nonlinear systems in the presence of noisy offline data and neglected nonlinearities. Assuming feasibility of \eqref{eqn:P2}, it is of interest to provide an estimate of the region of attraction of the resulting closed-loop system. We do this by following a similar procedure as in Section~\ref{sec:global_vs_local} (above Example~\ref{ex:ex_nominal}). In particular, suppose \eqref{eqn:P2} is feasible and consider the Lyapunov function $V(x_t)=x_t^\top P(x_t) x_t$ along with its difference
\begin{equation}
    \label{eqn:Vdiff_robust}
V(x_{t+1}) - V(x_t) = x_{t+1}^\top P(x_{t+1}) x_{t+1} - x_t^\top P(x_t) x_t.
\end{equation}
Notice from \eqref{eqn:robust_CL_NLsys_data} that
\[
    x_{t+1} = Y_1G(x_t)Z(x_t) + \underbrace{\left(\varepsilon^f(x_t) - D_0 G(x_t) Z(x_t)\right)}_{\eqqcolon r(x_t)},
\]
where $r(x)$ is unknown but bounded by a known function, i.e., $\norm{r(x)}\leq\bar{r}(x)$ for all $x\in\mathcal{X}\cap\mathcal{S}$, with
\[
\bar{r}(x) \coloneqq \bar\varepsilon + \sqrt{\lambda_{\max}(\Delta)}\norm{G(x)Z(x)},
\]
and $\lambda_{\max}$ denoting the maximum eigenvalue. Moreover, notice that $P(x_{t+1}) = \sum_{i=1}^N\alpha_i(x_{t+1})Q_i^{-1}$ and 
\[
\alpha_i(x_{t+1}) = \alpha_i(Y_1G(x_t)Z(x_t) + r(x_t)).
\]
Since $r(x)$ is unknown, one cannot evaluate $P(x_{t+1})$. However, due to the fact that $\alpha_i(x_{t+1})\geq0$ and $\sum_{i=1}^{N}\alpha_i(x_{t+1})=~1$, we can bound $P(x_{t+1})$ as follows
\begin{equation}
	P(x_{t+1}) =\sum\nolimits_{i=1}^N \alpha_i(x_{t+1})Q_i^{-1} \preceq \sum\nolimits_{i=1}^N Q_i^{-1}\eqqcolon \overline{P}.
	\label{eqn:Pbar}
\end{equation}
Plugging everything back into \eqref{eqn:Vdiff_robust} results in
\begin{align}
    &V(x_{t+1}) - V(x_t)\notag\\
    &\leq \left(Y_1G(x_t)Z(x_t) + r(x_t)\right)^\top \overline{P} \left(Y_1G(x_t)Z(x_t) + r(x_t)\right)\notag\\
    &\quad - x_t^\top P(x_t) x_t\notag\\
    &= \left(Y_1G(x_t)Z(x_t)\right)^\top \overline{P} \left(Y_1G(x_t)Z(x_t)\right) + r(x_t)^\top\overline{P}r(x_t)\notag\\
    &\quad + 2r(x_t)^\top \overline{P}Y_1G(x_t)Z(x_t) - x_t^\top P(x_t)x_t\notag\\
    &\leq \left(Y_1G(x_t)Z(x_t)\right)^\top \overline{P} \left(Y_1G(x_t)Z(x_t)\right) + \bar{r}^2(x_t)\norm{\overline{P}}\notag\\
    &\quad + 2\overline{r}(x_t)\norm{\overline{P}Y_1G(x_t)Z(x_t)} - x_t^\top P(x_t)x_t\notag\\
    &\eqqcolon \bar{h}(x_t), \label{eqn:robust_Vdiffbnd}
\end{align}
where $\bar{h}(x_t)$ is a known (and potentially conservative) upper bound that is composed of quantities which can be inferred from data and the solution of \eqref{eqn:P2}. As a result, any sub level set $\mathcal{R}_{c}=\left\lbrace x~|~ V(x)\leq c\right\rbrace$ contained in $\overline{\mathcal{V}}\cap\mathcal{X}\cap\mathcal{S}$, where $\overline{\mathcal{V}}\coloneqq~\left\lbrace x ~|~ \bar{h}(x) < 0 \right\rbrace\cup\{0\}$, is a positively invariant set and represents an estimate of the region of attraction of the closed-loop system, compare \cite{DePersis2023,haddad2008nonlinear}.

\begin{remark}
	A similar result to Corollary~\ref{cor:Kx} can be developed for the robust case considered in this section (not reported for space reasons). In particular, one can reduce the number of LMIs in \eqref{eqn:P2} from $N^2+N$ to $2N$ by using a uniform $Q$ matrix for all vertices. In this case, the above steps for computing an estimate of the ROA are simplified since $P(x_{t+1}) = P(x_t) = Q^{-1}$ as opposed to the upper bound computed in \eqref{eqn:Pbar}. However, this does not necessarily lead to a less conservative bound $\bar{h}(x)$ on the difference of the Lyapunov function, as it depends on the corresponding solution $Q$.\hfill$\square$
\end{remark}

\begin{figure}
    \centering
    \includegraphics[width=0.9\linewidth]{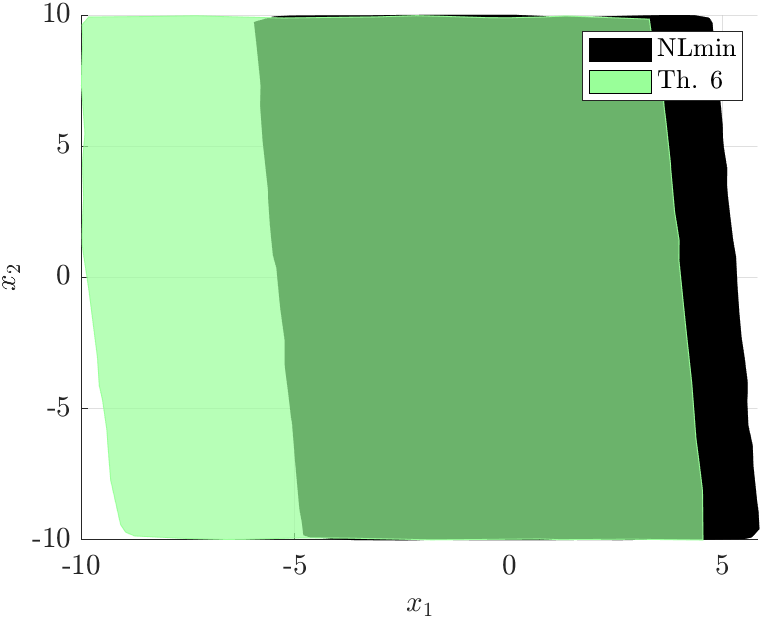}
    \caption{Numerical estimates of the true ROA of the closed-loop systems using (i) a polytopic controller obtained using Theorem~\ref{thm:Kx_robust} (in green) and (ii) a controller obtained from \cite{DePersis2023} (in black). The areas of each region are as follows: 213.6048 for the nonlinearity minimization approach and 272.7167 for the proposed approach corresponding to an increase of approximately 27.7\% in the size of the ROA.}
    \label{fig:robustROA}
\end{figure}

In the following, we illustrate the results of Theorem~\ref{thm:Kx_robust} with a numerical example and compare it to the nonlinearity minimization approach from \cite{DePersis2023}.
\begin{example}\label{ex:ex_robust}
    We now revisit the inverted pendulum system from Example~\ref{ex:global_inv_pend}, but instead consider the following choice of basis functions $Z(x)=\begin{bmatrix}x_{1} & x_{2} & x_1^2 & x_1^3\end{bmatrix}^\top$. In this case, Assumption~\ref{asmp:f2} is satisfied with $\bar\varepsilon = 0.005$ on the following subset of the state-space $\mathcal{S}\coloneqq\{ x ~|~ |x_1|\leq \pi/3\}$. For the above choice of basis functions, there exists $L$ of the form
    \[
    L(x_1) = \begin{bmatrix}
        x_1 & 0 & -1 & 0\\
        0 & 0 & x_1 & -1
    \end{bmatrix},
    \]
    which satisfies Assumption~\ref{asmp:L} (compare Theorem~\ref{thm:asmpLmonomoialZ} and Example~\ref{ex:ex_nominal}). We use $N=2$ vertices located at $\mathcal{L}_1 = L(0.5)$ and $\mathcal{L}_2= L(-0.5)$, respectively, and hence, Assumption (A2.3) is satisfied on $\mathcal{X}\coloneqq\{x\in\mathbb{R}^2~|~ |x_1|\leq 0.5\}$. We run an open-loop experiment of length $T=20$ and collect data by applying an input sampled from a uniform distribution $[-0.2,0.2]$ and initial conditions sampled from the same interval. The measured state is affected by uniform bounded noise with $\bar w = 0.001$ and, for a choice of $\Delta = 10^{-4} I$, it holds that $D_0D_0^\top\preceq \Delta$, where $D_0$ is defined below \eqref{eqn:noisy_state_identity}. 
    
    Using the same data set and the same parameter $\Delta$, we solve for (i) a polytopic controller using \eqref{eqn:P2} (where we set $\delta=0.01$) and (ii) a controller designed using the method from \cite[Eq. (49)]{DePersis2023} (there we set the user-defined parameter $\Omega=10^{-3}I$). Both programs are feasible and return stabilizing controllers. However, due to the conservative upper bound $\bar{h}(x)$ on the Lyapunov function difference in \eqref{eqn:robust_Vdiffbnd}, the procedure described above yields an empty ROA estimate (the same is true for the method from \cite{DePersis2023}). Instead, we numerically estimate the true ROA by running 5000 closed-loop experiments each of which is 200 steps long and record the initial conditions for which the experiment has converged. Convergence here is determined by testing whether the values of the state reached machine precision, i.e., $\leq 10^{-6}$ in 200 steps. Figure~\ref{fig:robustROA} depicts the set of all such initial conditions which empirically estimates the true ROA of the closed-loop system. It can be seen that the proposed framework provides a larger ROA estimate than the nonlinearity minimization approach from \cite{DePersis2023}.\hfill$\square$
\end{example}

\subsection{More general class of systems}
\label{sec:inaff}
So far, we only considered nonlinear systems of the form \eqref{eqn:NLsys}, which involve state-independent input vector fields (i.e., a constant $B$ matrix). In this subsection, we consider a more general class of nonlinear systems of the form 
\begin{equation}
	x_{t+1} =  f(x_t) + g(x_t)u_t,\label{eqn:generalNLsys}
\end{equation}
where $f:\mathbb{R}^n\to\mathbb{R}^n$ and $g:\mathbb{R}^n\to\mathbb{R}^{n\times m}$ are unknown continuous functions. The function $f$ satisfies Assumption~\ref{asmp:f}, while $g$ satisfies the following assumption.
\begin{assumption}\label{asmp:g}
	The function $g$ can be written as $g(x)=BW(x)$ where $B\in\mathbb{R}^{n\times q}$ is an unknown matrix of coefficients which satisfies $BB^\top\preceq \overline{B}$ for some known $\overline{B}\succeq0$, and $W:\mathbb{R}^{n}\to\mathbb{R}^{q\times m}$ is a known continuous matrix-valued function satisfying $W(x) = \sum_{i=1}^N\alpha_i(x) \mathcal{W}_i$ where $\mathcal{W}_i=W(v_i)$, $v_i\in\mathcal{X}$ and $\alpha_i$ as defined in Assumption~\ref{asmp:L}.\hfill$\square$
\end{assumption}

Such an assumption is satisfied, e.g., for affine functions $W$ with $\mathcal{X}$ being a polytope (compare \cite[p. 384]{toth2017handbook}). Alternatively, one can relax this assumption by considering $g(x)=BW(x)+\varepsilon^{g}(x)$ and then proceeding in a manner similar to Section~\ref{sec:inexact}. This, however, will not be presented here.

Assumptions~\ref{asmp:f} and \ref{asmp:g} allow us to write the system in the form $x_{t+1} = AZ(x_t)+BW(x_t)u_t$. Unlike previous sections, we now seek a controller of the form $u=KZ(x)$ such that the origin of the closed-loop system
\begin{equation}
	x_{t+1} = (A+BW(x_t)K)Z(x_t)\label{eqn:general_mbCLNLsys}
\end{equation}
is asymptotically stable. A polytopic gain matrix $K(x)$ as in \eqref{eqn:CLgain} is not applicable here since, in general, the product $W(x)K(x)$ destroys the affine polytopic structure of the closed–loop system (i.e., the dependence on $\alpha_i$ becomes quadratic). In this case, the admissibility conditions in Theorem~\ref{thm:LPV_polyquadadmissibility} (see \cite{bara2011dilated}) would no longer be applicable for the corresponding embedded LPVd system (compare \eqref{eqn:general_embeddedCL} below).

\begin{remark}
	Notice that a control law with a constant gain matrix is a special case of \eqref{eqn:CLgain} where all $K_i$'s are equal. In all previous sections, designing a control law with a constant gain matrix $K$ for systems of the form \eqref{eqn:NLsys} is possible by enforcing the decision variables $R_i$ in, e.g., \eqref{eqn:P1}, to be equal to one another. This, however, may lead to a smaller set of feasible solutions and worse closed-loop performance.\hfill$\square$
\end{remark}

As usual, we proceed by first deriving a data-based representation of the closed-loop system \eqref{eqn:general_mbCLNLsys}. To that end, we perform an experiment and collect input-state data $\{u_k^d,x_k^d\}_{k=0}^{T}$ that are then arranged in the following matrices
\begin{equation}
	\begin{aligned}
		&U_0 = \begin{bmatrix}
			u_0^d & u_1^d & \cdots & u_{T-1}^d
		\end{bmatrix}\in\mathbb{R}^{m\times T},\\
		&Z_0 = \begin{bmatrix}
			Z(x_0^d) & Z(x_1^d) & \cdots & Z(x_{T-1}^d)
		\end{bmatrix}\in\mathbb{R}^{S\times T},\\
		&\resizebox{\columnwidth}{!}{$\overline{W}_0 = \begin{bmatrix}
				W(x_0^d)u_0^d & W(x_1^d)u_1^d & \cdots & W(x_{T-1}^d)u_{T-1}^d
			\end{bmatrix}\in\mathbb{R}^{q\times T},$}\\
		&X_1 = \begin{bmatrix}
			x_1^d & x_2^d & \cdots & x_{T}^d
		\end{bmatrix}\in\mathbb{R}^{n\times T},  
	\end{aligned}
\end{equation}
which satisfy the following equation 
\begin{equation}
	X_1 = AZ_0 + B\overline{W}_0.\label{eqn:general_dataidentity}
\end{equation}
Assuming that the following holds
\begin{equation}
	\begin{bmatrix}
		U_0\\ Z_0
	\end{bmatrix}G =\begin{bmatrix}
	K\\ I_S
\end{bmatrix},\label{eqn:general_KG}
\end{equation}
the closed-loop system can be written as
\begin{equation}
	\begin{aligned}
		x_{t+1}&= (A+BW(x_t)K)Z(x_t)\\
		&\stackrel{\eqref{eqn:general_KG}}{=} (AZ_0G+BW(x_t)U_0G)Z(x_t)\\
		&\stackrel{\eqref{eqn:general_dataidentity}}{=} (X_1 - B\overline{W}_0 + BW(x_t)U_0)GZ(x_t)\\
		&\eqqcolon (X_1 + B\mathcal{F}(x_t))GZ(x_t)
	\end{aligned}\label{eqn:general_CL}
\end{equation}
where we have defined $\mathcal{F}(x_t)\coloneqq W(x_t)U_0-\overline{W}_0$. Notice that, since $W$ satisfies Assumption~\ref{asmp:g} and $\sum_{i=1}^{N}\alpha_i(x)=1$, then $\mathcal{F}(x_t)=\sum_{i=1}^{N}\alpha_i(x_t)\mathcal{F}_i$ with $\mathcal{F}_i\coloneqq\mathcal{W}_iU_0-\overline{W}_0$.

Analogously to Theorem~\ref{thm:embedding}, it can be shown that poly-quadratic admissibility of a polytopic LPVd system of the form
\begin{equation}
	\begin{bmatrix}
		I_n&0\\0 &0
	\end{bmatrix}\eta_{t+1} = \begin{bmatrix}
	(X_1 + B\mathcal{F}(\rho_t)) G\\ L(\rho_t)
\end{bmatrix}\eta_t,\label{eqn:general_embeddedCL}
\end{equation}
for which $\rho_t\in\mathcal{X}\subseteq\mathbb{R}^n$, implies asymptotic stability of the origin of \eqref{eqn:general_CL}. The following theorem represents an extension of Theorem~\ref{thm:Kx} to general input-affine systems and provides data-dependent LMI conditions for finding a stabilizing gain matrix $K$. For brevity, the proof is omitted as it follows similar steps as in the proof of\footnote{Observe that \eqref{eqn:general_CL} is similar to \eqref{eqn:robust_linearsubsystem} with $-D_0$ replaced by $+B\mathcal{F}(x)$. Therefore, one can follow the same steps as in the proof of Theorem~\ref{thm:Kx_robust}, with the exception that Petersen's lemma is applied here to eliminate the unknown $B$ instead of $D_0$.} Theorems~\ref{thm:Kx} and \ref{thm:Kx_robust}.
\begin{theorem}\label{thm:Kx_general}
	Consider a nonlinear system as in \eqref{eqn:generalNLsys} and let Assumptions~\ref{asmp:f},\ref{asmp:L} and \ref{asmp:g} hold. Fix some scalar $\delta\in\mathbb{R}$ and consider \eqref{eqn:Pgeneral} (shown at the top of the next page) in the decision variables $\lambda_i>0,\,R\in\mathbb{R}^{T\times S},\,Q_i=Q_i^\top\succ0\in\mathbb{R}^{n\times n}$, for all $i\in\mathcal{I}=\{1,\ldots,N\}$, $M\in\mathbb{R}^{n\times n},\,C_1\in\mathbb{R}^{S-n \times n}$ and an invertible $C_2\in\mathbb{R}^{S-n \times S-n}$.
	\begin{figure*}[!t]
		\centering
		\begin{subequations}
			\begin{align}
				&\begin{bmatrix}
					\sym\left(\begin{bmatrix}
						X_1\\ \mathcal{L}_i Z_0
					\end{bmatrix}R\begin{bmatrix}
						\delta I_{n} & 0\\0 & I_{S-n}
					\end{bmatrix}\right) + \begin{bmatrix}
						Q_j - \lambda_i \overline{B} & 0\\0 & 0
					\end{bmatrix} & \star & \star\\ 
					\left(\begin{bmatrix}
						X_1\\ \mathcal{L}_i Z_0
					\end{bmatrix}R\begin{bmatrix}
						I_{n}\\ 0
					\end{bmatrix}+\begin{bmatrix}
						\delta M^\top\\ 0
					\end{bmatrix}\right)^\top & \sym(M)-Q_i & \star\\
					-\mathcal{F}_i R\begin{bmatrix}
						\delta I_{n} & 0\\ 0 & I_{S-n}
					\end{bmatrix} & -\mathcal{F}_i R\begin{bmatrix}
						I_n\\0
					\end{bmatrix} & \lambda_i I_{q}
				\end{bmatrix}\succ0,\quad\forall (i,j)\in\mathcal{I}\times\mathcal{I},\label{eqn:Pgeneral_a}\\
				&Z_0 R = \begin{bmatrix}
					M & 0 \\ C_1 & C_2
				\end{bmatrix}.\label{eqn:Pgeneral_b}
			\end{align}\label{eqn:Pgeneral}
		\end{subequations}
		\hrule
	\end{figure*}
	If \eqref{eqn:Pgeneral} is feasible, then the origin of \eqref{eqn:general_CL}, with
	\begin{equation}
		K = U_0 R \begin{bmatrix}
			M & 0 \\ C_1 & C_2
		\end{bmatrix}^{-1},\label{eqn:Pgeneral_gain}
	\end{equation}
	is locally asymptotically stable. A Lyapunov function for the closed-loop system is given by $V(x_t)=x_t^\top P(x_t)x_t$ with $P(x_t) = \sum_{i=1}^{N}\alpha_i(x_t)Q_i^{-1}$. Furthermore, if $\mathcal{X}=\mathbb{R}^n$, then the origin is globally asymptotically stable.\hfill$\square$
\end{theorem}

As in previous sections, a similar result to Corollary~\ref{cor:Kx} can be derived here, i.e., by considering a uniform matrix $Q$ for all vertices. Furthermore, when $\mathcal{X}\subset\mathbb{R}^n$, data-driven estimates of the region of attraction can be obtained as follows. Consider the Lyapunov function $V$ and its difference along the system trajectory
\begin{equation}
	\begin{aligned}
		&V(x_{t+1}) - V(x_t) \\
		&= x_{t+1}^\top P(x_{t+1}) x_{t+1} - x_{t}^\top P(x_t) x_t\\
		&\resizebox{\columnwidth}{!}{$=\left((X_1 + B\mathcal{F}(x_t))GZ(x_t)\right)^\top P(x_{t+1})(X_1 + B\mathcal{F}(x_t))GZ(x_t)$}\\
		&\quad - x_{t}^\top P(x_t) x_t
	\end{aligned}\label{eqn:general_LyapDiff}
\end{equation}
Since $B$ is unknown, both $x_{t+1}$ and $P(x_{t+1})$ contain unknown terms, and hence \eqref{eqn:general_LyapDiff} cannot be evaluated directly. Thus, we bound the latter as in \eqref{eqn:Pbar} and make use of the fact that ${B}{B}^\top\preceq\overline{B}$ from Assumption~\ref{asmp:g} to write
\begin{equation}
	\begin{aligned}
		&V(x_{t+1}) - V(x_t) \\
		&\leq Z(x_t)^\top G^\top \left( 2 X_1^\top \overline{P} X_1 + c_1\mathcal{F}(x_t)^\top\mathcal{F}(x_t) \right) G Z(x_t)\\
		&\quad - x_t^\top P(x_t) x_t\\
		&\coloneqq \overline{h}(x_t)
	\end{aligned}
\end{equation}
where $c_1 = 2\lambda_{\max}(\overline{P})\lambda_{\max}(\overline{B})$. Here, $\overline{h}(x_t)$ is a known upper bound that is composed of quantities which can be inferred from data and the solution of \eqref{eqn:Pgeneral}. As a result, any sub level set $\mathcal{R}_{c}=\left\lbrace x~|~ V(x)\leq c\right\rbrace$ contained in $\overline{\mathcal{V}}\cap\mathcal{X}$, where $\overline{\mathcal{V}}\coloneqq\left\lbrace x ~|~ \overline{h}(x) < 0 \right\rbrace\cup\{0\}$, is a positively invariant set and represents an estimate of the region of attraction of the closed-loop system, compare \cite{DePersis2023,haddad2008nonlinear}.

\begin{figure}[!t]
	\centering
	\includegraphics[width=0.9\linewidth]{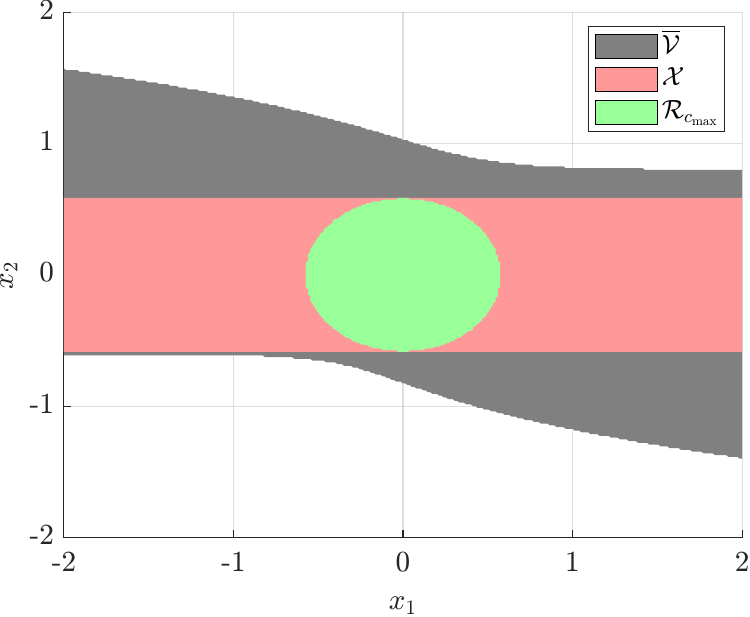}
	\includegraphics[width=0.95\linewidth]{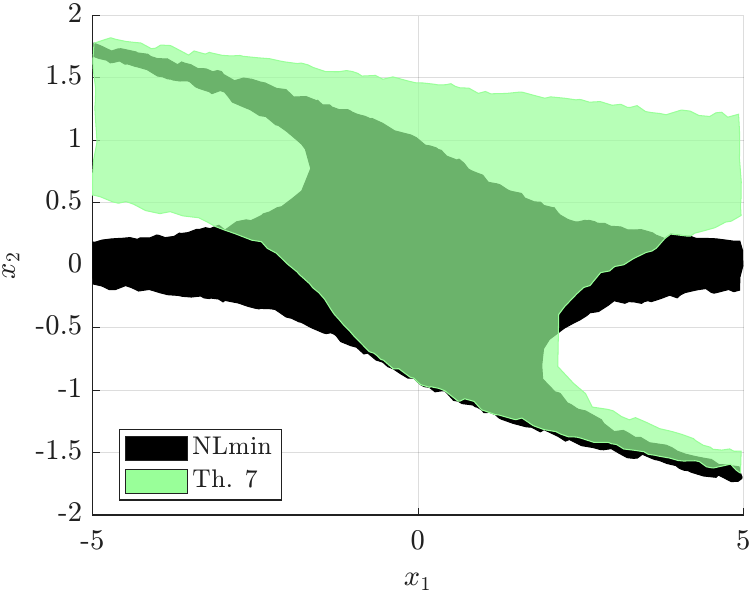}
	\caption{(Top): Data-dependent ROA estimate (in green) for the closed-loop system whose control law is designed using Theorem~\ref{thm:Kx_general}. In contrast, the method from \cite{guo2023data} returns an empty estimate of the ROA.\\ (Bottom): Numerical estimates of the true ROA of the closed-loop systems using (i) a controller obtained using Theorem~\ref{thm:Kx_general} (in green) and (ii) a controller obtained from \cite{guo2023data} (in black). The areas of each region are as follows: 11.8030 for the nonlinearity minimization approach and 17.6380 for the proposed approach corresponding to an increase of approximately 49.4\% in the size of the ROA.}
	\label{fig:generalsystems}
\end{figure}

In the following we illustrate the results of Theorem~\ref{thm:Kx_general} with an example and compare it to that of \cite{guo2023data}, where the nonlinearity minimization approach from \cite{DePersis2023} was extended to the class of general input-affine systems.
\begin{example}
	We consider the following system (previously considered in \cite{guo2023data})
	\[
	\begin{aligned}
		x_{1,t+1} &= 0.5 x_{2,t},\\
		x_{2,t+1} &= x_{1,t} + x_{2,t}^3 + (1+x_{2,t})u_t.
	\end{aligned}
	\]
	We use $Z(x)=\begin{bmatrix}
		x^\top & x_2^2 & x_2^3
	\end{bmatrix}^\top$ and $W(x)=\begin{bmatrix}
	1 & x_2
\end{bmatrix}^\top$. Analogously to Theorem~\ref{thm:asmpLmonomoialZ}, one can construct $L$ of the form
\[
L(x_2) = \begin{bmatrix}
	0 & x_2 & -1 & 0\\
	0 & 0 & x_2 & -1
\end{bmatrix}.
\]
All assumptions are satisfied on the set $\mathcal{X} = \{x\in\mathbb{R}^2 ~|~ |x_2|\leq 0.6\}\subset\mathbb{R}^2$. We run an open-loop experiment of length $T=10$ and collect data by applying an input sampled from a uniform distribution $[-0.3,0.3]$ and initial conditions sampled from the same interval. We set $\delta=0$ and find that \eqref{eqn:Pgeneral} is feasible. Using the same data set, we follow the procedure in \cite{guo2023data} to obtain a controller of the form $u=K_{\mathrm{NLmin}}Z(x)$. To compare the performance of the two controllers, we estimate the region of attraction of the closed-loop systems using data only. The results are reported in Figure~\ref{fig:generalsystems} (see next page). We point out that the method from \cite{guo2023data} returns an empty estimate of the ROA when using data only. This may be attributed to the fact that the effect of the nonlinearities in the state-dependent input vector fields are simply minimized, which necessitates conservative upper bounds on the Lyapunov function difference to estimate the ROA. In contrast, our approach exploits the inherent nonlinearities of the system and returns a non-empty estimate of the ROA (following the procedure described below Theorem~\ref{thm:Kx_general}). 

We additionally compare numerical estimates of the true ROA obtained by running 5000 closed-loop experiments each of which is 200 steps long and record the initial conditions for which the experiment has converged. Convergence here is determined by testing whether the values of the state reached machine precision, i.e., $\leq 10^{-6}$ in 200 steps. Figure~\ref{fig:generalsystems} additionally depicts the set of all such initial conditions which empirically estimates the true ROA of the closed-loop system. It can be seen that the proposed framework provides a larger ROA estimate than the approach from \cite{guo2023data}.\hfill$\square$
\end{example}
\section{Discussion and conclusions}
\label{sec:conc}
In this paper, we introduced a notion of descriptor embedding for discrete-time nonlinear systems. Such a notion opens the possibilities to analyze and design controllers for nonlinear systems, both in model-based and data-based domains, by utilizing the rich theory of descriptor systems. We illustrate in this paper how one can use admissibility results of discrete-time linear parameter-varying descriptor systems to (globally) asymptotically stabilize the origin of nonlinear systems. Specifically, we provide sufficient data-dependent LMI conditions that, if feasible, return a stabilizing controller. 

Notably, our method allows for obtaining globally stabilizing controllers for nonlinear systems without resorting to nonlinearity cancellation. This is in contrast to existing techniques in, e.g., \cite{DePersis2023,guo2023data} which, on the other hand, offer more interpretability and less computational complexity. Additionally, we extend our results to account for neglected nonlinearities and noisy data. When only local stabilization can be guaranteed, we provide estimates for the corresponding region of attraction using data only. 

Our results show that the proposed method is competitive with existing techniques in the literature. In cases when cancellation is not possible or in the presence of measurement noise and inexact basis function expansion, the proposed method was shown to outperform existing techniques. This is mainly due to the proposed design procedure which makes use of the available nonlinearities instead of minimizing and/or dominating their effect.

The notion of descriptor embedding has the potential for many future research directions. For instance, it is of interest to investigate optimal control techniques for nonlinear systems from the point of view of optimal control of descriptor systems, compare \cite{Bender87}. Apart from set-point stabilization, other interesting problems to investigate include output regulation as well as robust invariance or tracking problems.
\appendix
\subsection{On PDLF structure for poly-quadratically admissible LPVd systems}\label{app:onPDLF}
In \cite{bara2011dilated}, sufficient LMI conditions were provided to test for poly-quadratic admissibility of polytopic LPVd systems (see Theorem~\ref{thm:LPV_polyquadadmissibility}). It is of interest for us to obtain a PDLF for the corresponding subsystem \eqref{eqn:bareta1}. More specifically, we wish to show that \eqref{eqn:PDLFinequalitybareta1} holds with $\bar{P}(\rho_t)=\sum_{i=1}^{N}\alpha_i(\rho_t)Q_i^{-1}$. To this end, we revisit the proof of \cite[Thm. 4.3(ii)]{bara2011dilated}. There it is shown that \eqref{eqn:polyquadraticadmissibilityLMIconditions} implies
\[\begin{bmatrix}
	\sym\left(\bar{A}(\rho_t)F\right) + Q(\rho_{t+1}) & \bar{A}(\rho_t)G + F^\top \\
	G^\top\bar{A}(\rho_t)^\top + F & \sym(G) - Q(\rho_t)
\end{bmatrix}\succ0,\]
where $Q(\rho_{t+1}) = \sum_{j=1}^{N}\alpha_j(\rho_{t+1}) Q_j$ and $Q(\rho_t) = \sum_{i=1}^{N}\alpha_i(\rho_{t}) Q_i$. Pre- and post- multiplying by $\begin{bmatrix} I_r & - \bar{A}(\rho_t)\end{bmatrix}$ and its transpose results in
\[
Q(\rho_{t+1}) - \bar{A}(\rho_t)Q(\rho_t)\bar{A}(\rho_{t})^\top\succ0.
\]
Since the above holds for every convex combination (recall the above definition of $Q(\rho_{t+1}), Q(\rho_t)$), then it also holds at the vertices, i.e.,
\[Q_j - \bar{A}(\rho_t) Q_i \bar{A}(\rho_t)^\top \succ0, \quad \forall(i,j)\in\mathcal{I}\times\mathcal{I}.\]
Using the Schur complement twice, it holds that
\[
Q_i^{-1} - \bar{A}(\rho_t)^\top Q_j^{-1} \bar{A}(\rho_t) \succ 0,\quad \forall(i,j)\in\mathcal{I}\times\mathcal{I}.
\]
Letting $\bar{P}_i \coloneqq Q_i^{-1}$, $\bar{P}_j \coloneqq Q_j^{-1}$, multiplying by $\alpha_i(\rho_t)\alpha_j(\rho_{t+1})$ and summing over $i,j\in\mathcal{I}$ results in
\[\bar{P}(\rho_t) - \bar{A}(\rho_t)^\top \bar{P}(\rho_{t+1}) \bar{A}(\rho_t) \succ 0,\]
where we have defined $\bar{P}(\rho_t) \coloneqq \sum_{i=1}^{N}\alpha_{i}(\rho_t) \bar{P}_i$, $\bar{P}(\rho_{t+1}) \coloneqq \sum_{j=1}^{N}\alpha_{j}(\rho_{t+1}) \bar{P}_j$ and made use of the fact that $\alpha_i,\alpha_j$ sum to one. This last inequality is precisely \eqref{eqn:PDLFinequalitybareta1}, as desired.

\subsection{Example of a function \texorpdfstring{$\mathnormal{L}$}{L} satisfying Assumption~\ref{asmp:L}}\label{app:exmpL}
Consider a choice of $Z$ composed of all monomials of a two-dimensional state $x=[x_1 \,\, x_2]^\top$ up to degree three. Then following the procedure described in the proof of Theorem~\ref{thm:asmpLmonomoialZ} we obtain
\[
\resizebox{\linewidth}{!}{$\underbrace{\begin{bmatrix}
    x_1 & 0 & -1 & 0 & 0 & 0 & 0 & 0 & 0\\
    0 & x_2 & 0 & -1 & 0 & 0 & 0 & 0 & 0\\
    x_2 & 0 & 0 & 0 & -1 & 0 & 0 & 0 & 0\\
    0 & 0 & x_1 & 0 & 0 & -1 & 0 & 0 & 0\\
    0 & 0 & 0 & x_2 & 0 & 0 & -1 & 0 & 0\\
    0 & 0 & 0 & 0 & x_1 & 0 & 0 & -1 & 0\\
    0 & 0 & 0 & 0 & x_2 & 0 & 0 & 0 & -1
\end{bmatrix}}_{=L(x)}\begin{bmatrix}
    x_1 \\ x_2 \\ x_1^2 \\ x_2^2 \\ x_1x_2 \\ x_1^3 \\ x_2^3 \\ x_1^2x_2 \\ x_1x_2^2
\end{bmatrix} = 0.$}
\]
Here, we have
\[
\begin{aligned}
L(x) &= \left[\begin{array}{c|c}L_1(x) & L_2(x)\end{array}\right]\\
&=\left[
\begin{array}{cc|ccccccc}
x_1 & 0   & -1 & 0  & 0  & 0  & 0  & 0 & 0\\
0   & x_2 & 0  & -1 & 0  & 0  & 0  & 0 & 0\\
x_2 & 0   & 0  & 0  & -1 & 0  & 0  & 0 & 0\\
0   & 0   & x_1& 0  & 0  & -1 & 0  & 0 & 0\\
0   & 0   & 0  & x_2& 0  & 0  & -1 & 0 & 0\\
0   & 0   & 0  & 0  & x_1& 0  & 0  & -1& 0\\
0   & 0   & 0  & 0  & x_2& 0  & 0  & 0 & -1
\end{array}
\right],
\end{aligned}
\]
where $L_2(x)$ is invertible for all $x$. Finally, since $L$ is an affine function, its image over a polytopic domain $\mathcal{X}$ can be expressed as (see \cite[p. 384]{toth2017handbook})
\[
L(x) = \sum\nolimits_{i=1}^N \alpha_i(x) \mathcal{L}_i, \quad \alpha_i(x) \geq 0,\quad\sum\nolimits_{i=1}^N \alpha_i(x) = 1,
\]
where $\mathcal{L}_i=L(v_i),\,v_i\in\mathcal{X}$. The functions $\alpha_i$ can be taken as the generalized barycentric coordinates (see \cite{floater2015generalized}) with respect to the polytope $\mathcal{X}$. For this example, let $\mathcal{X} = \lbrace x ~|~ |x_i| \leq \gamma_i \rbrace\subset\mathbb{R}^2$ whose $N=4$ vertices given by
\[
\begin{aligned}
	\mathcal{L}_1 &= L([\gamma_1 \,\, \gamma_2]^\top),\\
	\mathcal{L}_3 &= L([\gamma_1 \,\, -\gamma_2]^\top),
\end{aligned}\quad
\begin{aligned}
	\mathcal{L}_2 &= L([-\gamma_1 \,\, \gamma_2]^\top),\\
	\mathcal{L}_4 &= L([-\gamma_1 \,\, -\gamma_2]^\top).
\end{aligned}
\]
One can then choose $\alpha_i$, e.g., as 
\[
\begin{aligned}
    \alpha_1(x) &= \frac{(\gamma_1+x_1)(\gamma_2+x_2)}{4\gamma_1\gamma_2},\\
    \alpha_3(x) &= \frac{(\gamma_1+x_1)(\gamma_2-x_2)}{4\gamma_1\gamma_2},
\end{aligned}\,\,
\begin{aligned}
    \alpha_2(x) &= \frac{(\gamma_1-x_1)(\gamma_2+x_2)}{4\gamma_1\gamma_2},\\
    \alpha_4(x) &= \frac{(\gamma_1-x_1)(\gamma_2-x_2)}{4\gamma_1\gamma_2}.
\end{aligned}
\]

\section*{References}
\bibliographystyle{IEEEtran}
\bibliography{references}

\begin{thebibliography}{10}
\providecommand{\url}[1]{#1}
\csname url@samestyle\endcsname
\providecommand{\newblock}{\relax}
\providecommand{\bibinfo}[2]{#2}
\providecommand{\BIBentrySTDinterwordspacing}{\spaceskip=0pt\relax}
\providecommand{\BIBentryALTinterwordstretchfactor}{4}
\providecommand{\BIBentryALTinterwordspacing}{\spaceskip=\fontdimen2\font plus
\BIBentryALTinterwordstretchfactor\fontdimen3\font minus
  \fontdimen4\font\relax}
\providecommand{\BIBforeignlanguage}[2]{{%
\expandafter\ifx\csname l@#1\endcsname\relax
\typeout{** WARNING: IEEEtran.bst: No hyphenation pattern has been}%
\typeout{** loaded for the language `#1'. Using the pattern for}%
\typeout{** the default language instead.}%
\else
\language=\csname l@#1\endcsname
\fi
#2}}
\providecommand{\BIBdecl}{\relax}
\BIBdecl

\bibitem{markovsky2021behavioral}
I.~Markovsky and F.~D{\"o}rfler, ``Behavioral systems theory in data-driven
  analysis, signal processing, and control,'' \emph{Annual Reviews in Control},
  vol.~52, pp. 42--64, 2021.

\bibitem{florian23ddcontrol}
F.~Dörfler, ``Data-driven control: Parts {I} and {II},'' \emph{IEEE Control
  Systems Magazine}, vol.~43, no. 5, 6, pp. 24--27, 27--31, 2023.

\bibitem{willems1997book}
J.~C. Willems and J.~W. Polderman, \emph{Introduction to mathematical systems
  theory: a behavioral approach}.\hskip 1em plus 0.5em minus 0.4em\relax
  Springer Science \& Business Media, 1997, vol.~26.

\bibitem{willems2005note}
J.~C. Willems, P.~Rapisarda, I.~Markovsky, and B.~L. De~Moor, ``A note on
  persistency of excitation,'' \emph{Systems \& Control Letters}, vol.~54,
  no.~4, pp. 325--329, 2005.

\bibitem{markovsky2022identifiability}
I.~Markovsky and F.~D{\"o}rfler, ``Identifiability in the behavioral setting,''
  \emph{IEEE Transactions on Automatic Control}, vol.~68, no.~3, pp.
  1667--1677, 2022.

\bibitem{Alsalti2025DBrep}
M.~Alsalti, I.~Markovsky, V.~G. Lopez, and M.~A. Müller, ``Data-based system
  representations from irregularly measured data,'' \emph{IEEE Transactions on
  Automatic Control}, vol.~70, no.~1, pp. 143--158, 2025.

\bibitem{depersis2019formulas}
C.~De~Persis and P.~Tesi, ``Formulas for data-driven control: Stabilization,
  optimality, and robustness,'' \emph{IEEE Transactions on Automatic Control},
  vol.~65, no.~3, pp. 909--924, 2019.

\bibitem{vanwaarde2023informativity}
H.~J. Van~Waarde, J.~Eising, M.~K. Camlibel, and H.~L. Trentelman, ``The
  informativity approach: To data-driven analysis and control,'' \emph{IEEE
  Control Systems Magazine}, vol.~43, no.~6, pp. 32--66, 2023.

\bibitem{bisoffi2022data}
A.~Bisoffi, C.~De~Persis, and P.~Tesi, ``Data-driven control via {P}etersen’s
  lemma,'' \emph{Automatica}, vol. 145, p. 110537, 2022.

\bibitem{van2020noisy}
H.~J. van Waarde, M.~K. Camlibel, and M.~Mesbahi, ``From noisy data to feedback
  controllers: Nonconservative design via a matrix {S}-lemma,'' \emph{IEEE
  Transactions on Automatic Control}, vol.~67, no.~1, pp. 162--175, 2020.

\bibitem{berberich2024overview}
J.~Berberich and F.~Allgöwer, ``An overview of systems-theoretic guarantees in
  data-driven model predictive control,'' \emph{Annual Review of Control,
  Robotics, and Autonomous Systems}, vol.~8, no. Volume 8, 2025, pp. 77--100,
  2025.

\bibitem{alsalti2023notes}
M.~Alsalti, V.~G. Lopez, and M.~A. Müller, ``Notes on data-driven
  output-feedback control of linear {MIMO} systems,'' \emph{IEEE Transactions
  on Automatic Control}, vol.~70, no.~9, pp. 6143--6150, 2025.

\bibitem{Li2024}
L.~Li, A.~Bisoffi, C.~{De Persis}, and N.~Monshizadeh, ``Controller synthesis
  from noisy-input noisy-output data,'' \emph{Automatica}, vol. 183, p. 112545,
  2026.

\bibitem{berberich2020trajectory}
J.~Berberich and F.~Allg{\"o}wer, ``A trajectory-based framework for
  data-driven system analysis and control,'' in \emph{2020 European Control
  Conference (ECC)}, 2020, pp. 1365--1370.

\bibitem{rueda2020data}
J.~G. Rueda-Escobedo and J.~Schiffer, ``Data-driven internal model control of
  second-order discrete {V}olterra systems,'' in \emph{2020 59th IEEE
  Conference on Decision and Control (CDC)}, 2020, pp. 4572--4579.

\bibitem{verhoek2021fundamental}
C.~Verhoek, R.~T{\'o}th, S.~Haesaert, and A.~Koch, ``Fundamental lemma for
  data-driven analysis of linear parameter-varying systems,'' in \emph{2021
  60th IEEE conference on decision and control (CDC)}, 2021, pp. 5040--5046.

\bibitem{alsalti2023data}
M.~Alsalti, V.~G. Lopez, J.~Berberich, F.~Allg{\"o}wer, and M.~A. M{\"u}ller,
  ``Data-based control of feedback linearizable systems,'' \emph{IEEE
  Transactions on Automatic Control}, vol.~68, no.~11, pp. 7014--7021, 2023.

\bibitem{markovsky2023bilinear}
I.~Markovsky, ``Data-driven simulation of generalized bilinear systems via
  linear time-invariant embedding,'' \emph{IEEE Transactions on Automatic
  Control}, vol.~68, no.~2, pp. 1101--1106, 2023.

\bibitem{martin2023guarantees}
T.~Martin, T.~B. Sch{\"o}n, and F.~Allg{\"o}wer, ``Guarantees for data-driven
  control of nonlinear systems using semidefinite programming: A survey,''
  \emph{Annual Reviews in Control}, p. 100911, 2023.

\bibitem{depersis2023annrev}
C.~De~Persis and P.~Tesi, ``Learning controllers for nonlinear systems from
  data,'' \emph{Annual Reviews in Control}, p. 100915, 2023.

\bibitem{Guo2022polynomial}
M.~Guo, C.~De~Persis, and P.~Tesi, ``Data-driven stabilization of nonlinear
  polynomial systems with noisy data,'' \emph{IEEE Transactions on Automatic
  Control}, vol.~67, no.~8, pp. 4210--4217, 2022.

\bibitem{DePersis2023}
C.~De~Persis, M.~Rotulo, and P.~Tesi, ``Learning controllers from data via
  approximate nonlinearity cancellation,'' \emph{IEEE Transactions on Automatic
  Control}, vol.~68, no.~10, pp. 6082--6097, 2023.

\bibitem{hu2024enforcing}
Z.~Hu, C.~De~Persis, and P.~Tesi, ``Enforcing contraction via data,''
  \emph{IEEE Transactions on Automatic Control}, 2025, {E}arly Access.

\bibitem{mauroy2020}
Mauroy, I.~Mezi{\'c}, and Y.~Susuki, \emph{{K}oopman Operator in Systems and
  Control: Concepts, Methodologies, and Applications}.\hskip 1em plus 0.5em
  minus 0.4em\relax Springer International Publishing, 2020.

\bibitem{robin2024koopman}
R.~Strässer, M.~Schaller, K.~Worthmann, J.~Berberich, and F.~Allgöwer,
  ``Koopman-based feedback design with stability guarantees,'' \emph{IEEE
  Transactions on Automatic Control}, vol.~70, no.~1, pp. 355--370, 2025.

\bibitem{Haseli2022}
M.~Haseli and J.~Cortés, ``Learning {K}oopman eigenfunctions and invariant
  subspaces from data: Symmetric subspace decomposition,'' \emph{IEEE
  Transactions on Automatic Control}, vol.~67, no.~7, pp. 3442--3457, 2022.

\bibitem{dai1989}
L.~Dai, \emph{Singular control systems}.\hskip 1em plus 0.5em minus 0.4em\relax
  Springer, 1989.

\bibitem{verhoek2023general}
C.~Verhoek, P.~J.~W. Koelewijn, S.~Haesaert, and R.~Tóth, ``Direct data-driven
  state-feedback control of general nonlinear systems,'' in \emph{2023 62nd
  IEEE Conference on Decision and Control (CDC)}, 2023, pp. 3688--3693.

\bibitem{bara2011dilated}
G.~I. Bara, ``Dilated {LMI} conditions for time-varying polytopic descriptor
  systems: the discrete-time case,'' \emph{International Journal of Control},
  vol.~84, no.~6, pp. 1010--1023, 2011.

\bibitem{barbosa2018}
K.~A. Barbosa, C.~E. de~Souza, and D.~Coutinho, ``Admissibility analysis of
  discrete linear time-varying descriptor systems,'' \emph{Automatica},
  vol.~91, pp. 136--143, 2018.

\bibitem{daafouz2001parameter}
J.~Daafouz and J.~Bernussou, ``Parameter dependent {L}yapunov functions for
  discrete time systems with time varying parametric uncertainties,''
  \emph{Systems \& control letters}, vol.~43, no.~5, pp. 355--359, 2001.

\bibitem{de2008robust}
C.~E. De~Souza, K.~A. Barbosa, and M.~Fu, ``Robust filtering for uncertain
  linear discrete-time descriptor systems,'' \emph{Automatica}, vol.~44, no.~3,
  pp. 792--798, 2008.

\bibitem{Alsalti_pe}
M.~Alsalti, V.~G. Lopez, and M.~A. Müller, ``On the design of persistently
  exciting inputs for data-driven control of linear and nonlinear systems,''
  \emph{IEEE Control Systems Letters}, vol.~7, pp. 2629--2634, 2023.

\bibitem{haddad2008nonlinear}
W.~M. Haddad and V.~Chellaboina, \emph{Nonlinear dynamical systems and control:
  a Lyapunov-based approach}.\hskip 1em plus 0.5em minus 0.4em\relax Princeton
  university press, 2008.

\bibitem{toth2017handbook}
C.~D. Toth, J.~O'Rourke, and J.~E. Goodman, \emph{Handbook of discrete and
  computational geometry}.\hskip 1em plus 0.5em minus 0.4em\relax CRC press,
  2017.

\bibitem{Stone1948}
M.~H. Stone, ``The generalized {W}eierstrass approximation theorem,''
  \emph{Mathematics Magazine}, vol.~21, no.~4, pp. 167--184, 1948.

\bibitem{floater2015generalized}
M.~S. Floater, ``Generalized barycentric coordinates and applications,''
  \emph{Acta Numerica}, vol.~24, pp. 161--214, 2015.

\bibitem{masubuchi1997h}
I.~Masubuchi, Y.~Kamitane, A.~Ohara, and N.~Suda, ``${H}_{\infty}$ control for
  descriptor systems: A matrix inequalities approach,'' \emph{Automatica},
  vol.~33, no.~4, pp. 669--673, 1997.

\bibitem{Boyd2004}
S.~P. Boyd and L.~Vandenberghe, \emph{Convex optimization}.\hskip 1em plus
  0.5em minus 0.4em\relax Cambridge university press, 2004.

\bibitem{petersen1987}
I.~R. Petersen, ``A stabilization algorithm for a class of uncertain linear
  systems,'' \emph{Systems \& control letters}, vol.~8, no.~4, pp. 351--357,
  1987.

\bibitem{guo2023data}
M.~Guo, C.~De~Persis, and P.~Tesi, ``Data-driven control of input-affine
  systems via approximate nonlinearity cancellation,''
  \emph{IFAC-PapersOnLine}, vol.~56, no.~2, pp. 1357--1362, 2023.

\bibitem{Bender87}
D.~J. Bender and A.~J. Laub, ``The linear-quadratic optimal regulator for
  descriptor systems: Discrete-time case,'' \emph{Automatica}, vol.~23, no.~1,
  pp. 71--85, 1987.

\end{thebibliography}

\begin{IEEEbiography}[{\includegraphics[width=1in,height=1.25in,clip,keepaspectratio]{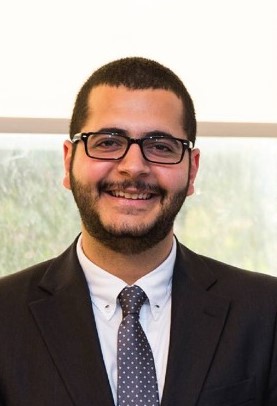}}]{Mohammad Alsalti} received his B.Sc. in Mechanical Engineering from the University of Jordan, Jordan, in 2017. In 2017, he was an intern at NASA Ames Research Center as part of the intelligent robotics group. In 2018, Mohammad was awarded the Fulbright scholarship to pursue graduate studies. In 2020, he obtained his M.Sc. in Mechanical Engineering from the University of Maryland, College Park, USA. He is currently a Ph.D. student at the Institute of Automatic Control at Leibniz University Hannover, Germany. He is working on developing data-driven control techniques for linear and nonlinear systems. In 2024, he was a visiting researcher at Prof. Claudio De Persis's research group at the University of Groningen, the Netherlands.
\end{IEEEbiography}

\begin{IEEEbiography}[{\includegraphics[width=1in,height=1.25in,clip,keepaspectratio]{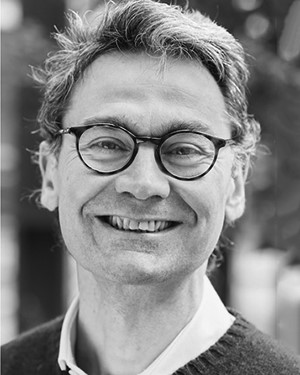}}]{Claudio De Persis} received the Laurea and Ph.D. degrees in engineering from the University of Rome ``La Sapienza'', Rome, Italy, in 1996 and 2000, respectively. 
	
	He is currently a Professor with the Engineering and Technology Institute, University of Groningen, Groningen, the Netherlands, since 2011. From 2000 to 2001, he held Postdoctoral positions with Washington University, St. Louis, MO, USA, from 2001 to 2002, with Yale University, New Haven, CT, USA, and from 2002 to 2009 and 2009 to 2011, he held faculty positions with the University of Rome ``La Sapienza'', and Twente University, Enschede, the Netherlands. His main research interests include automatic control and its applications.
\end{IEEEbiography}

\begin{IEEEbiography}[{\includegraphics[width=1in,height=1.25in,clip,keepaspectratio]{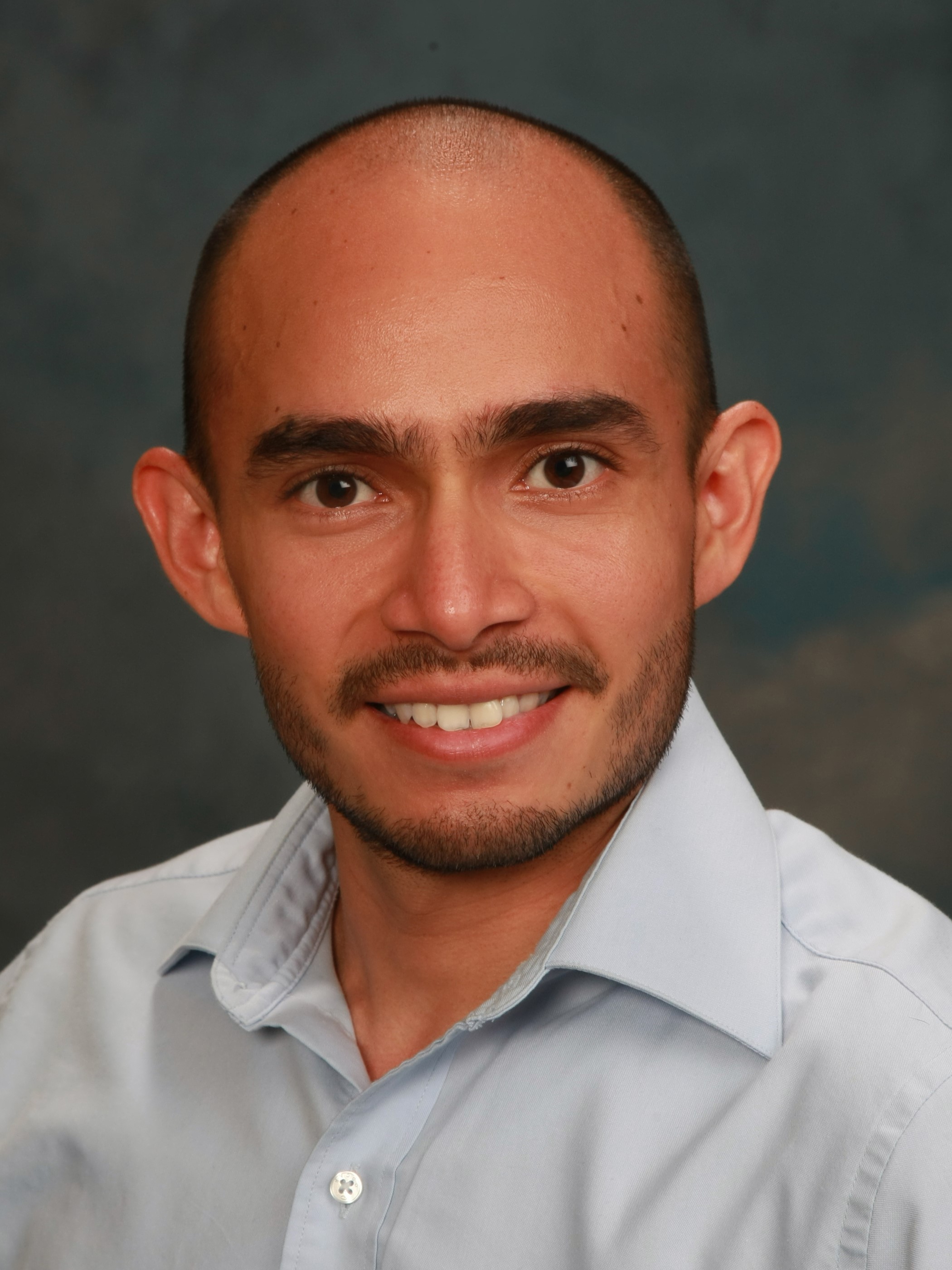}}]{Victor G. Lopez}
	received his B.Sc. degree in Communications and Electronics Engineering from the Universidad Autonoma de Campeche, in Campeche, Mexico, in 2010, the M.Sc. degree in Electrical Engineering from the Research and Advanced Studies Center (Cinvestav), in Guadalajara, Mexico, in 2013, and his Ph.D. degree in Electrical Engineering from the University of Texas at Arlington, Texas, USA, in 2019. In 2015 Victor was a Lecturer at the Western Technological Institute of Superior Studies (ITESO) in Guadalajara, Mexico. From August 2019 to June 2020, he was a postdoctoral researcher at the University of Texas at Arlington Research Institute and an Adjunct Professor in the Electrical Engineering department at UTA. Victor is currently a postdoctoral researcher at the Institute of Automatic Control, Leibniz University Hannover, in Hannover, Germany. His research interest include cyber-physical systems, reinforcement learning, game theory, distributed control and robust control. 
\end{IEEEbiography}

\begin{IEEEbiography}[{\includegraphics[width=1in,height=1.25in,clip,keepaspectratio]{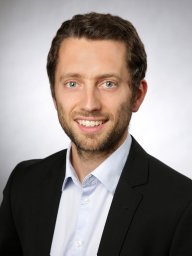}}]{Matthias A. Müller} received a Diploma degree in engineering cybernetics from the University of Stuttgart, Germany, an M.Sc. in electrical and computer engineering from the University of Illinois at Urbana-Champaign, US (both in 2009), and a Ph.D. from the University of Stuttgart in 2014. Since 2019, he is Director of the Institute of Automatic Control and Full Professor at the Leibniz University Hannover, Germany. 
	
His research interests include nonlinear control and estimation, model predictive control, and data- and learning-based control, with application in different fields including biomedical engineering and robotics. He has received various awards for his work, including the 2015 European Systems \& Control PhD Thesis Award, the inaugural Brockett-Willems Outstanding Paper Award for the best paper published in Systems \& Control Letters in the period 2014-2018, an ERC starting grant in 2020, the IEEE CSS George S. Axelby Outstanding Paper Award 2022, and the Journal of Process Control Paper Award 2023. He serves/d as an associate editor for Automatica and as an editor of the International Journal of Robust and Nonlinear Control.
\end{IEEEbiography}

\end{document}